\newcommand{\M}{{\mathcal M}} 
\newcommand{\R}{I\!\!R}   
\newcommand{\N}{I\!\!N}   
\newcommand{\Z}{{\mathbb Z}} 
  
\newcommand{\sgn}{\mbox{\sl sgn}}
\newcommand{\dgn}{\mbox{\sl dgn}}
\newcommand{\rk}{\mbox{\sl rk}}

\newcommand{\Bsym}{{\mathrm B}_{\mathrm{sym}}}
\newcommand{\gr}{g^{(\mathrm r)}}
\newcommand{\maslov}{\mathrm i_{\scriptscriptstyle\rm M}} 
\newcommand{\corpo}{{I\!\!K}}
\newcommand{\poincare}{\mathfrak P}
\newcommand{\rip}[2]{\langle\!\!\langle #1,#2\rangle\!\!\rangle}
\newcommand{\ddt}{\frac{\mathrm d}{\mathrm dt}}
\newcommand{\tto}{\big\vert_{t=t_0}}

\documentclass[oneside,draft,a4paper,11pt]{amsart}

\usepackage{amsmath}               
\usepackage{amsthm}                
\usepackage{times}

\title[A Generalized Index Theorem in semi-Riemannian Geometry]%
{A Generalized Index Theorem for Morse-Sturm Systems 
and Applications to semi-Riemannian Geometry}
\author[F.\ Giannoni]{Fabio Giannoni}
\address{Dipartimento di Matematica e Fisica,
Universit\'a di Camerino, Italy}
\email{giannoni@campus.unicam.it}
\author[A.\ Masiello]{Antonio Masiello}
\address{Dipartimento Interuniversitario di Matematica,
Politecnico di Bari, Italy}
\email{masiello@pascal.dm.uniba.it}
\author[P.\ Piccione]{Paolo Piccione}
\address{Departamento de Matem\'atica,  Universidade de S\~ao Paulo, Brazil}
\email{piccione@ime.usp.br}
\urladdr{http://www.ime.usp.br/\~{}piccione}
\author[D.\ Tausk]{Daniel V.\ Tausk}
\address{Departamento de Matem\'atica,  Universidade de S\~ao Paulo, Brazil}
\email{tausk@ime.usp.br}

\begin{document}


\theoremstyle{definition}\newtheorem{defsym}{Definition}[section]
\theoremstyle{definition}\newtheorem{separating}[defsym]{Definition} 
\theoremstyle{plain}\newtheorem{C1}[defsym]{Lemma}
\theoremstyle{definition}\newtheorem{defC1subspaces}[defsym]{Definition}
\theoremstyle{plain}\newtheorem{HSelementary}[defsym]{Proposition}
\theoremstyle{remark}\newtheorem{remdegenerate}[defsym]{Remark}
\theoremstyle{remark}\newtheorem{remconstant}[defsym]{Remark}
\theoremstyle{plain}\newtheorem{corelementary}[defsym]{Corollary}
\theoremstyle{plain}\newtheorem{produce}[defsym]{Lemma}

\theoremstyle{plain}\newtheorem{discrete}{Proposition}[section]
\theoremstyle{remark}\newtheorem{remId+K}[discrete]{Remark}
\theoremstyle{plain}\newtheorem{BC1}[discrete]{Proposition}
\theoremstyle{plain}\newtheorem{positive}[discrete]{Proposition}
\theoremstyle{plain}\newtheorem{corpositive}[discrete]{Corollary}
\theoremstyle{plain}\newtheorem{central}[discrete]{Proposition}
\theoremstyle{remark}\newtheorem{remdiscrete}[discrete]{Remark}
\theoremstyle{plain}\newtheorem{RiemannianMorse}[discrete]{Corollary}

\theoremstyle{definition}\newtheorem{defadmissible}{Definition}[section]
\theoremstyle{plain}\newtheorem{indexthm}[defadmissible]{Theorem}
\theoremstyle{plain}\newtheorem{FtC1}[defadmissible]{Lemma}
\theoremstyle{plain}\newtheorem{surjective}[defadmissible]{Lemma}
\theoremstyle{plain}\newtheorem{corC1}[defadmissible]{Corollary}
\theoremstyle{plain}\newtheorem{corsoma}[defadmissible]{Corollary}
\theoremstyle{plain}\newtheorem{orthogonal}[defadmissible]{Lemma}
\theoremstyle{plain}\newtheorem{corKer}[defadmissible]{Corollary}
\theoremstyle{plain}\newtheorem{Idcomp}[defadmissible]{Proposition}
\theoremstyle{plain}\newtheorem{initial}[defadmissible]{Proposition}
\theoremstyle{plain}\newtheorem{summarize}[defadmissible]{Theorem}

\theoremstyle{plain}\newtheorem{extension}{Theorem}[section]

\theoremstyle{plain}\newtheorem{MorseLorentzian}{Theorem}[section]
\theoremstyle{remark}\newtheorem{exsimple}[MorseLorentzian]{Example}
\theoremstyle{remark}\newtheorem{excausal}[MorseLorentzian]{Example}
\theoremstyle{remark}\newtheorem{exparallel}[MorseLorentzian]{Example}
\theoremstyle{remark}\newtheorem{remsemi}[MorseLorentzian]{Remark}
\theoremstyle{plain}\newtheorem{dim2}[MorseLorentzian]{Proposition}

\theoremstyle{definition}\newtheorem{defcprec}{Definition}[section]
\theoremstyle{plain}\newtheorem{globalMorse}[defcprec]{Theorem}


\begin{abstract} 
We prove an extension of the Index Theorem for Morse--Sturm systems
of the form $-V''+RV=0$, where $R$ is symmetric with respect to a (non positive) 
symmetric bilinear form, and thus the corresponding differential operator 
is not self-adjoint.
The result is then applied to the case of a Jacobi equation along
a geodesic in a Lorentzian manifold, obtaining an extension of
the Morse Index Theorem for Lorentzian geodesics with variable
initial endpoints. 
Given a Lorentzian manifold $(\M,g)$, we consider a geodesic $\gamma$ in $\M$
starting orthogonally to a smooth
submanifold $\mathcal P$ of $\M$. Under suitable hypotheses, satisfied, 
for instance, if $(\M,g)$ is stationary, the
theorem  gives an equality between the index of the second variation of the 
action functional $f$ at $\gamma$ and the sum of the {\em Maslov index\/} of 
$\gamma$ with the index of the metric $g$ on $\mathcal P$. 
Under generic circumstances, the Maslov index of $\gamma$ is given by an algebraic 
count of the $\mathcal P$-focal points along~$\gamma$.
Using the Maslov index, we obtain the global Morse relations for geodesics
between two fixed points in a stationary Lorentzian manifold.
\end{abstract}

\maketitle

\begin{section}{Introduction}
\label{sec:intro}
The goal of this paper is to prove an index theorem for Morse--Sturm systems
of differential equations with coefficients that are symmetric with respect 
to an indefinite inner product of $\R^n$. The main motivation for this
kind of investigation comes from semi-Riemannian geometry, where Morse--Sturm
systems appear in the form of {\em Jacobi\/} equations for vector fields
along geodesics.
 
Let $(\M,g)$ be a semi-Riemannian manifold, $\mathcal P$ a smooth
submanifold of $\M$ and $\gamma:[0,1]\mapsto\M$ be
a geodesic in $\M$, with $\gamma(0)\in\mathcal P$ and
$\dot\gamma(0)\in T_{\gamma(0)}\mathcal P^\perp$; set $q=\gamma(1)$. 
The curve $\gamma$ is then a stationary point of the {\em action\/} functional 
\[f(z)=\frac12\int_0^1g(\dot z,\dot z)\;\mathrm dt,\] defined in the space $\Omega_{\mathcal P,q}$ of  
curves joining $\mathcal P$ and the point $q$ in $\mathcal M$.
The {\em index form\/} $I_{\{\gamma,\mathcal P\}}$ is the symmetric bilinear
form  given by the second variation of $f$,
defined on the tangent space $T_\gamma \Omega_{\mathcal P,q}$,
which consists of vector fields $V$ along  $\gamma$ with $V(0)\in T_{\gamma(0)}\mathcal P$
and $V(1)=0$. We recall the definition of $I_{\{\gamma,\mathcal P\}}$:
\begin{equation}\label{eq:indexform}
\begin{split}
I_{\{\gamma,\mathcal P\}}(V,W)=&\int_0^1\Big[g(\nabla_{\dot\gamma}V,\nabla_{\dot\gamma}W)
+g(R(\dot\gamma,V)\,\dot\gamma,W)\Big]\;\mathrm dt+\\&-S_{\dot\gamma(0)}(V(0),W(0)),
\end{split}
\end{equation}
where $\nabla$ is the covariant derivative of the Levi--Civita connection of $g$,
$R$ is the curvature tensor of $\nabla$ and $S_{\dot\gamma(0)}$ is the second
fundamental form of $\mathcal P$ in the direction of $\dot\gamma(0)$.

One obtains an infinite dimensional Hilbertian structure
in $\Omega_{\mathcal P,q}$ by requiring a Sobolev $H^1$-regularity for the curves
in $\Omega_{\mathcal P,q}$; then, $I_{\{\gamma,\mathcal P\}}$ is a {\em bounded\/}
bilinear symmetric form on the Hilbert space $T_\gamma \Omega_{\mathcal P,q}$.

If $(\M,g)$ is Riemannian, i.e., if $g$ is a positive
definite metric tensor, the celebrated Morse Index Theorem (see for instance
\cite[Theorem~2.2]{dC}, \cite[Theorem~15.1]{M}, \cite{PT} ) 
states that the {\em index\/} of $I_{\{\gamma,\mathcal P\}}$, which is the
dimension of a maximal subspace of $T_\gamma \Omega_{\mathcal P,q}$ on which 
$I_{\{\gamma,\mathcal P\}}$ is negative definite, equals
the {\em geometric index\/} $\mathrm i_{\mathrm{geom}}(\gamma)$
of $\gamma$, which is the number of
$\mathcal P$-focal points along $\gamma$ counted with multiplicity. 
Such equality can also be given in terms of the multiplicity of the
negative eigenvalues of the Jacobi differential operator, which 
is a self-adjoint operator representing the index form in the Hilbert
space of square-integrable vector fields along $\gamma$.

From the viewpoint of Calculus
of Variations, the elements of $T_\gamma \Omega_{\mathcal P,q}$ are 
interpreted as {\em infinitesimal variations\/}
of $\gamma$, and the index of $I_{\{\gamma,\mathcal P\}}$ on 
$T_\gamma \Omega_{\mathcal P,q}$ is {\em the number of
essentially different directions in which $\gamma$ can be deformed in order to obtain a curve 
of shorter length}.

The theorem has been successively extended by Beem and Ehrlich
to Lorentzian manifolds (see \cite{BE, BEE}),
i.e., manifolds endowed with a metric tensor $g$ of index $1$, in the
case of {\em causal\/} (non spacelike) geodesics. For such an extension
one only needs minor modifications to the original statement (and proof) of
the theorem. Most notably one needs to consider the restriction
of $I_{\{\gamma,\mathcal P\}}$ to the space $T_\gamma \Omega_{\mathcal P,q}^\perp$ 
of vector fields along $\gamma$ which
are pointwise orthogonal to $\gamma$. With this restriction, which in
the Riemannian case is totally ininfluent for the computation of the
index of $I_{\{\gamma,\mathcal P\}}$, one basically excludes the variations of $\gamma$
obtained by simple {\em reparameterizations\/} of $\gamma$. For timelike
Lorentzian geodesics, the affine parameterization is the one that
{\em maximizes\/} the value of the action functional, and thus 
the restriction to $T_\gamma \Omega_{\mathcal P,q}^\perp$ has the effect of factoring out from $T_\gamma
\Omega_{\mathcal P,q}$ an infinite dimensional space on which $I_{\{\gamma,\mathcal P\}}$ is
negative definite, thus making the restricted $I_{\{\gamma,\mathcal P\}}$ into a form with finite index.

For spacelike Lorentzian geodesics, or more in general for geodesic
of any causal character in semi-Riemannian manifolds with metrics
of index greater than or equal to two, there is no hope to extend
the original formulation of the index theorem, due mainly to the 
following reasons:
\begin{itemize}
\item the index of $I_{\{\gamma,\mathcal P\}}$ on both $T_\gamma \Omega_{\mathcal P,q}$ and 
$T_\gamma \Omega_{\mathcal
P,q}^\perp$ is infinite;
\item the set of $\mathcal P$-focal points along a geodesic may fail  to be
discrete, and there is no meaningful  notion of geometric index;
\item the Jacobi differential operator is no longer self-adjoint.
\end{itemize} 
In the case of a geodesic $\gamma$ having only a finite number of $\mathcal P$-focal
points, one can ask the question of whether there exists a {\em natural\/}
subspace $\mathcal K^\gamma$ of $T_\gamma \Omega_{\mathcal P,q}$ with the property that the 
restriction of $I_{\{\gamma,\mathcal P\}}$
to $\mathcal K^\gamma$ has finite index, equal to the geometric index of $\gamma$.
However, also for this special case the question seems to have
a negative answer, due to the fact that, while the index of a bilinear
form has some (semi-)continuity properties, the geometric index
is {\em not stable\/} by small perturbations. Indeed, one
can produce examples where (isolated) $\mathcal P$-focal points simply {\em evaporate\/}
by arbitrary small perturbations of the metric (see~\cite{MPT}), or
examples of a sequence $\gamma_n$ of geodesics having a finite number
of $\mathcal P$-focal points converging to a geodesic $\gamma$ that has a {\em continuum\/}
of $\mathcal P$-focal points (see \cite{Hel}).

In order to prove an extension of the index theorem
in semi-Riemannian geometry one needs to determine a natural subspace $\mathcal K^\gamma$
of the Hilbert space $T_\gamma\Omega_{\mathcal P,q}$ with the properties that:
\begin{itemize}
\item the index of the restriction of $I_{\{\gamma,\mathcal P\}}$ to $\mathcal K^\gamma$ is finite;
\item such index should be related to some geometrical properties
of the geodesic $\gamma$ and of the manifold $\mathcal P$.
\end{itemize}
A hint for the choice of such a space was given by recent studies (see \cite{GP,Ma})
concerning the geodesical connectedness of Lorentzian manifolds $(\M,g)$ whose metric $g$ is
{\em stationary\/}, i.e., there exists a globally defined Killing timelike vector field
on $\M$. Given any such vector field $Y$ on $\M$, one has a conservation law for
geodesics given by:
\begin{equation}\label{eq:conservation}
g(Y,\dot\gamma)\equiv c_\gamma\ \text{(constant)}.
\end{equation}
Considering the Hilbertian structure on $\Omega_{\mathcal P,q}$,
one proves that the set $\Omega_{\mathcal P,q}^Y$  of curves in
$\Omega_{\mathcal P,q}$  satisfying \eqref{eq:conservation} almost everywhere
is a smooth submanifold of $\Omega_{\mathcal P,q}$, and
that the critical points of the restriction of the action functional $f$ to 
$\Omega_{\mathcal P,q}^Y$ are precisely the geodesics joining $\mathcal P$ and $q$
in $\mathcal M$.
Given one such geodesic $\gamma$, the tangent space $\mathcal K^\gamma=
T_\gamma\Omega_{\mathcal P,q}^Y$ is the Hilbert subspace of $T_\gamma\Omega_{\mathcal P,q}$
consisting of those vector fields $V$ along $\gamma$ that satisfy the linearization
of \eqref{eq:conservation}. Using the Killing property of $Y$, 
the space $\mathcal K^\gamma$  can be described as:
\begin{equation}\label{eq:defNgamma}
\mathcal K^\gamma=\Big\{V\in T_\gamma\Omega_{\mathcal
P,q}:g(\nabla_{\dot\gamma}V,Y)-g(V,\nabla_{\dot\gamma}Y)\equiv C_V\ 
\text{(constant)}\Big\}.
\end{equation}
Using compact embeddings of the Sobolev space $H^1$ into the space
$C^0$, one then proves that the restriction of the index form $I_{\{\gamma,\mathcal P\}}$
to $\mathcal K^\gamma$ is represented by a self-adjoint operator, which is
a compact perturbation of the identity. In particular, its index is finite.
The definition of the space $\mathcal K^\gamma$ makes perfectly sense
also in the case that $Y$ is a timelike Jacobi field along $\gamma$, and
also in this case we have finiteness of the index of the restriction of
$I_{\{\gamma,\mathcal P\}}$ to $\mathcal K^\gamma$.
Observe that the restriction of a Killing field along a geodesic is Jacobi, 
and thus this second construction is more general.
This construction gives a solution for the first point mentioned in the
program above; the next step is to give a geometrical interpretation
of the value of the index of $I_{\{\gamma,\mathcal P\}}$ on $\mathcal K^\gamma$.

Inspired by some techniques in Hamiltonian systems (see \cite{Ar}), it has
recently been defined the notion of {\em Maslov index\/} for
a semi-Riemannian geodesic (see~\cite{Hel} and also \cite{MPT}),
which is an integer number given by a certain topological invariant.
Under {\em generic\/} circumstances, the Maslov index can be
computed as a sort of {\em algebraic count\/} of the multiplicities
of the $\mathcal P$-focal points. In particular, for Riemannian and causal
Lorentzian geodesics it is always equal to the geometric index
(see~\cite{MPT}).  
For spacelike Lorentzian geodesics, or more in general for all kinds
of geodesics in semi-Riemannian manifolds with metric tensor of
index greater or equal to two, the contribution of each $\mathcal P$-focal
point to the value of the index is an integer number, possibly zero
or negative, called the {\em signature\/} of the $\mathcal P$-focal point, 
whose absolute value is less than or equal to the multiplicity
of the $\mathcal P$-focal point. Generically, the Maslov index of a semi-Riemannian
geodesic is the sum of the signatures of its $\mathcal P$-focal points, and
this sum is in absolute value less than or equal to the geometric index
of the geodesic.
Besides the finiteness, a remarkable property
of the Maslov index is its {\em stability\/} by small perturbations
(see~\cite{MPT}), due to its topological nature.

In this sense, the Maslov index of a geodesic is a natural candidate for 
substituting the notion of geometric index for Riemannian and causal
Lorentzian geodesics.

The main result of the paper (Theorem~\ref{thm:extension} and its
geometrical formulation Theorem~\ref{thm:MorseLorentzian})
is that, if $\gamma(1)$ is not a 
$\mathcal P$-focal point along $\gamma$, then  the  index  of the
restriction of $I_{\{\gamma,\mathcal P\}}$ to $\mathcal K^\gamma$
is equal to the sum of the Maslov index of $\gamma$ and the index
of the restriction of the Lorentzian metric $g$ to $T_{\gamma(0)}\mathcal P$.
In particular, this number is independent on the choice of the vector field $Y$.
To strengthen the analogy with the classical index Theorem, we remark that
it was recently proven (see \cite[Theorem~6.2.3]{MPT}) that,
under generic circumstances,
the Maslov index of $\gamma$ is equal to the {\em spectral index\/} of 
$\gamma$, which is computed 
as a sort of algebraic count of the (real) negative eigenvalues of the Jacobi
differential operator.

When comparing with the classical result of the Morse index theorem
in Riemannian manifolds, we see that for non positive
definite metrics some  new phenomena appear:
\begin{itemize}
\item if $\mathcal P$ is timelike at $\gamma(0)$, i.e., if
the restriction of $g$ to $T_{\gamma(0)}\mathcal P$ has positive index, 
then the {\em initial value\/}
of the index of $I_{\{\gamma,\mathcal P\}}$ is strictly positive, hence even
 {\em small portions\/} of  $\gamma$ are never local minimizers for the 
restricted action  functional;
\item each $\mathcal P$-focal point along $\gamma$ gives a contribution to the index
which may be positive, negative or even null;
\item the multiplicity of the $\mathcal P$-focal points is not stable by
perturbations, and arbitrary small perturbations of a given geodesic
may create or destroy focal points (see~\cite{MPT}).
\end{itemize}
By a parallel trivialization of the tangent bundle of $\M$ along the geodesic
$\gamma$, one can reformulate the entire theory in terms of Morse--Sturm--Liouville
systems of differential equations in $\R^n$. In this framework, the version of the 
Index Theorem discussed in this paper may be considered an extension
of the Sturm Oscillation Theorem.

The proof of the main result of the paper is based on a general
method for computing the variation of the index of a smooth curve
$B(t)$ of symmetric bounded bilinear forms defined on a smooth
family $\mathcal H_t$ of Hilbert spaces (Proposition~\ref{thm:HSelementary}). 
The {\em jumps\/} of the index function $i(t)=\mathrm{ind}(B(t)\vert_{\mathcal H_t})$
occur at the instants where $B(t)$ becomes singular, that correspond to the
conjugate points. The value of the jump at a discontinuity point $t_0$ is then proven
to be equal to the signature of the corresponding conjugate point 
(Proposition~\ref{thm:central}), under the assumption that the derivative
$B'(t_0)$ be non degenerate on $\mathrm{Ker}(B(t_0))$.
Under these circumstances, such calculation gives the proof of the aimed index Theorem.

Finally, we need to emphasize the fact that the {\em stability\/} of the Morse index
and of the Maslov index (see~\cite{MPT}) plays a crucial role in the proof of our results. Namely,
in order to employ the method described, we need to make a technical assumption
concerning the non degeneracy of the restriction of $g$ to suitable subspaces.
Such assumption, which holds {\em generically}, is needed to guarantee the finiteness
of the set of conjugate points and it is the core of the proof of 
Proposition~\ref{thm:HSelementary}, where we show how to compute the  jump 
of the index function at each conjugate point. The proof of the general case
is then given using a {\em perturbation\/} argument, which is based on
the observation that both the Morse index and the Maslov index of a semi-Riemannian
geodesic do not change by small $C^0$-perturbations of the data.  

Some examples and applications of the theory developed are discussed 
in the final part of the paper. In particular, under a suitable completeness
assumption, we obtain the global Morse relations for geodesics with fixed endpoints 
in a stationary  Lorentzian manifold (Theorem~\ref{thm:globalMorse}).

For a {\em standard static\/} Lorentzian manifold, the Morse relations
have been proven in \cite{BM} using the Morse index of the energy
functional restricted to the set of curves satisfying the constraint
\eqref{eq:conservation}; the same kind of relations have been proven in \cite{GM}
in the more general case of a {\em standard stationary\/} metric in
a manifold with (possibly non smooth) convex boundary.

\end{section}
\begin{section}{Abstract Results in Functional Analysis}
\label{sec:abstract}
Given Banach spaces $E_1$ and $E_2$, we denote by $\mathcal L(E_1,E_2)$
the set of all bounded linear operators from $E_1$ to $E_2$ and by $\mathrm B(E_1,E_2,\R)$
the set of all bounded bilinear maps from $E_1\times E_2$ to $\R$. If $E_1=E_2=E$, we also
set $\mathcal L(E)=\mathcal L(E,E)$ and $\mathrm B(E,\R)=\mathrm B(E,E,\R)$; by $\mathrm
B_{\mathrm{sym}}(E,\R)$ we mean the set of symmetric bounded bilinear maps on $E$.

We give some general definitions 
concerning symmetric bilinear forms for later use.
\begin{defsym}\label{thm:defsym}
Let $V$ be any real vector space and $B:V\times V\mapsto\R$ a symmetric
bilinear form. The {\em negative type number\/} (or {\em index}) $n_-(B)$ of $B$
is the possibly infinite number defined by
\begin{equation}\label{eq:def-}
n_-(B)=\sup\Big\{{\rm dim}(W):W\ \text{subspace of}\ V \ \text{on which}\ B\ 
\text{is negative definite}\Big\}.
\end{equation}
The {\em positive type number\/} $n_+(B)$ is given by
$n_+(B)=n_-(-B)$; if at least one of these two numbers is finite,
the {\em signature\/} $\sgn(B)$ is defined by:
\[\sgn(B)=n_+(B)-n_-(B).\]
The {\em kernel\/} of $B$, ${\rm Ker}(B)$, is the set  of vectors $v\in V$
such that $B(v,w)=0$ for all $w\in V$; the {\em degeneracy\/} $\dgn(B)$
of $B$ is the (possibly infinite) dimension of ${\rm Ker}(B)$. 
\end{defsym}
If $V=V_+\oplus V_-$, where $B$ is positive semidefinite on $V_+$ and negative
definite on $V_-$, then $n_-(B)=\mathrm{dim}(V_-)$; for, obviously $n_-(B)\ge
\mathrm{dim}(V_-)$ and every subspace $S$ on which $B$ is negative definite
satisfies $S\cap V_+=\{0\}$, and therefore $\mathrm{dim}(S)\le\mathrm{dim}(V_-)$.
Moreover, if in addition $B$ is positive definite on $V_+$, then $\mathrm{Ker}(B)=\{0\}$.
Namely, if $v=v_++v_-\in\mathrm{Ker}(B)$, with $v_+\in V_+$ and $v_-\in V_-$, then, by considering
the equality $-B(v_+,v_-)=B(v_+,v_+)=B(v_-,v_-)$, we get $v_+=v_-=0$.
A simple density argument shows that if the symmetric bilinear form 
$B$ is continuous with respect to some norm in the vector space $V$, 
then its index does not change
when one extends $B$ to the Banach space completion of $V$.

If $V$ is finite dimensional, then the numbers $n_+(B)$, $n_-(B)$ and
$\dgn(B)$ are respectively  the number of $1$'s, $-1$'s and $0$'s
in the canonical form of $B$ as given by the Sylvester's Inertia Theorem.
In this case, $n_+(B)+n_-(B)$ is equal to the codimension of ${\rm Ker}(B)$,
and it is also called the {\em rank} of $B$, $\rk(B)$.
\smallskip

Given a Hilbert space $\mathcal H$ with inner product $\langle\cdot,\cdot\rangle$, 
to any bounded bilinear form $B:\mathcal H\times \mathcal H\mapsto\R$ by Riesz's theorem
there corresponds a bounded linear operator $T_B:\mathcal H\mapsto\mathcal H$, which
is related to $B$ by:
\begin{equation}\label{eq:defBT}
B(x,y)=\langle T_B(x),y\rangle,\quad\forall\,x,y\in\mathcal H.
\end{equation}
We say that $T_B$ is the linear operator {\em associated\/} to $B$ with respect
to the inner product $\langle\cdot,\cdot\rangle$. Clearly, $B$ is symmetric
if and only if $T_B$ is self-adjoint. We say that $B$ is {\em non degenerate\/}
if $T_B$ is injective; $B$ will be said to be {\em strongly non degenerate\/}
if $T_B$ is an isomorphism. If $T_B$ is a {\em Fredholm operator\/} of index
$0$, i.e., if $T_B$ is a compact perturbation of an isomorphism, then, by
the Fredholm's Alternative, $B$ is non degenerate if and only if it is strongly
non degenerate. Observe that the strong non degeneracy is stable by small
perturbations, since the set of isomorphisms of $\mathcal H$ is open
in $\mathcal L(\mathcal H)$.

\smallskip

We now give a criterion for the differentiability of curves
in Banach spaces. We start with a definition
\begin{separating}\label{thm:separating}
Let $E$ and $E_0$ be real Banach spaces. A subset $\Phi\subset\mathcal L(E,E_0)$
is said to be {\em separating\/} for $E$ if for all $x\in E\setminus\{0\}$ 
there exists $\phi\in\Phi$ such
that $\phi(x)\ne0$.
\end{separating}
We now prove the following:
\begin{C1}\label{thm:C1}
Let $E, E_0$ be real Banach spaces and $F,G:[a,b]\mapsto E$ be fixed maps, 
with $G$ continuous. Let $\Phi\subset\mathcal L(E,E_0)$ be a separating set 
for $E$; assume that for each $\phi\in\Phi$ the composition $\phi\circ F:[a,b]\mapsto E_0$
is of class $C^1$, and that $(\phi\circ F)'(t)=\phi\circ G(t)$ for all $t\in [a,b]$.
Then, $F$ is a map of class $C^1$, and $F'(t)=G(t)$ for all $t\in [a,b]$.
\end{C1}
\begin{proof}
Fix $t\in [a,b]$; we have to prove that $F'(t)=G(t)$.
We claim that the following equality holds:
\begin{equation}\label{eq:claim}
F(t+h)-F(t)=\int_t^{t+h}G(s)\,\mathrm ds.
\end{equation}
It follows easily by applying each element $\phi\in\Phi$ to both sides
of \eqref{eq:claim} and using the separating property of $\Phi$. Denoting by
$\Vert\cdot\Vert$ the norm of $E$, it follows:
\[\left\Vert \frac{F(t+h)-F(t)}h-G(t)\right\Vert\le \left\vert\frac1h
\int_t^{t+h}\Vert G(s)-G(t)\Vert\,\mathrm ds\right\vert;\]
the continuity of $G$ concludes the argument.
\end{proof}
In the next proposition and its corollary we exhibit a method to compute the variation
of the index of a curve of symmetric bilinear forms. We want to leave the domains
of the forms variable, and we use the following notion of a $C^1$-curve of closed subspaces
of a Hilbert space:
\begin{defC1subspaces}\label{thm:defC1subspaces}
Let $\mathcal H$ be a Hilbert space, $I\subset\R$ an interval and $\{\mathcal D_t\}_{t\in I}$
be a family of closed subspaces of $\mathcal H$. We say that $\{\mathcal D_t\}_{t\in I}$
is a {\em $C^1$-family\/} of subspaces if for all $t_0\in I$ there exists a $C^1$-curve
$\alpha:\,]t_0-\varepsilon,t_0+\varepsilon\,[\,\cap\, I\mapsto \mathcal L(\mathcal H)$ 
and a closed subspace $\overline{\mathcal D}\subset\mathcal H$
such that $\alpha(t)$ is an isomorphism  and $\alpha(t)(\mathcal D_t)=\overline{\mathcal D}$
for all~$t$.
\end{defC1subspaces}
We will call the maps $\alpha$ appearing in Definition~\ref{thm:defC1subspaces}
the {\em local trivializations\/} of the family $\{\mathcal D_t\}_{t\in I}$.

In the following Proposition we study how the index of a smooth curve
$B(t)$ of symmetric bilinear forms varies after passing through a degenerate
instant $t_0$. We need a technical assumption on the map $B(t_0)$, which
must be represented by a compact perturbation of a positive operator.
\begin{HSelementary}\label{thm:HSelementary}
Let $\mathcal H$ be a real Hilbert space with inner product
$\langle\cdot,\cdot\rangle$, and let $B:[t_0,t_0+r]\mapsto \mathrm
B_{\mathrm{sym}}(\mathcal H,\R)$, $r>0$, be a map of class $C^1$. 
Let $\{\mathcal D_t\}_{t\in[t_0,t_0+r]}$ be a $C^1$-family of closed subspaces
of $\mathcal H$, and denote by $\overline B(t)$ the restriction of $B(t)$
to $\mathcal D_t\times\mathcal D_t$.
Assume that the following three hypotheses are satisfied:
\begin{enumerate}
\item\label{itm:hp2.5.1} $\overline B(t_0)$ is represented by an operator
of the form $L+K$, with $L:\mathcal D_{t_0}\mapsto\mathcal D_{t_0}$ a positive
isomorphism and  $K:\mathcal D_{t_0}\mapsto\mathcal D_{t_0}$ a (self-adjoint) 
compact operator;
\item\label{itm:hp2.5.2} the restriction $\widetilde B$ of the derivative $B'(t_0)$ to 
$\mathrm{Ker}(\overline B(t_0))\times \mathrm{Ker}(\overline B(t_0))$ 
is non degenerate;
\item\label{itm:hp2.5.3} $\mathrm{Ker}(\overline B(t_0))\subseteq \mathrm{Ker}(B(t_0))$.
\end{enumerate}
Then, for
$t>t_0$ sufficiently close to $t_0$, $\overline B(t)$ is non degenerate, and we have:
\begin{equation}\label{eq:changen-}
n_-(\overline B(t))=n_-(\overline B(t_0))+n_-(\widetilde B),
\end{equation}
all the terms of the above equality being finite natural numbers.
\end{HSelementary}
\begin{proof}
By possibly passing to a smaller $r$, we can assume the existence of
a $C^1$-curve $\alpha(t)$ of isomorphisms of $\mathcal H$ such that $\alpha(t)$
carries $\mathcal D_t$ to a fixed subspace $\overline{\mathcal D}$ of $\mathcal H$.
We can now replace each $B(t)$ by the push-forward 
$B(t)(\alpha(t)^{-1}\cdot,\alpha(t)^{-1}\cdot)$,
and each $\mathcal D_t$ by $\overline{\mathcal D}$. Such replacements will not
affect the hypotheses of the Proposition, nor the quantities involved
 in the equality \eqref{eq:changen-}. For instance, thanks to the hypothesis
\ref{itm:hp2.5.3}, the index of the restriction of $B'(t_0)$ to
$\mathrm{Ker}(\overline B(t_0))$ does not change; namely, for $V,W\in\mathrm{Ker}
\Big(B(t_0)(\alpha(t_0)^{-1}\,\cdot\,,\alpha(t_0)^{-1}\,\cdot\,)\vert_{\overline{\mathcal
D}\times\overline{\mathcal D}}\Big)$, it is:
\begin{equation}\label{eq:accorti}
\begin{split}
&\ddt\,  B(t)(\alpha(t)^{-1}V,\alpha(t)^{-1}W)\tto =B'(t_0)(\alpha(t_0)^{-1}V,\alpha(t_0)^{-1}W)+\\
&+B(t_0)(\ddt\,\alpha(t)^{-1}V,\alpha(t_0)^{-1}W)\tto+B(t_0)(\alpha(t_0)^{-1}V,\ddt\alpha(t)^{-1}W)
\tto=\\ &\qquad\qquad\qquad=B'(t_0)(\alpha(t_0)^{-1}V,\alpha(t_0)^{-1}W).
\end{split}
\end{equation}
We can therefore assume without loss  of generality that $\mathcal D_t=\mathcal H$
and $\overline B(t)=B(t)$ for all $t$. Moreover, we observe here that, by a convenient
choice of the Hilbert space inner product on $\mathcal H$, we can assume that 
$\overline B(t_0)=B(t_0)$ is represented by a compact perturbation of
the {\em identity\/} of $\mathcal H$, $\mathrm{Id}+K$.

Now, the subspace $N=\mathrm{Ker}(B(t_0))$ is 
the eigenspace of $K$ corresponding to
the eigenvalue $-1$, hence it is finite dimensional. 

We start considering the case that $B(t_0)$ is positive semi-definite
on $\mathcal H$ and that $\widetilde B$ is positive definite on $N$.
In this case, the thesis means that $B(t)$ is positive definite
on $\mathcal H$ for $t>t_0$ sufficiently close to $t_0$.

Let $S$ be any closed complementary subspace of $N$ in $\mathcal H$; 
clearly $B(t_0)$ is positive definite on $S$. We claim that there exists
a positive constant $c_0$ such that, for $t$ sufficiently close to $t_0$,
it is:
\begin{equation}\label{eq:claim2}
\phantom{\quad\forall\,x\in S\ \text{with}\ \Vert x\Vert=1.}
B(t)[x,x]\ge c_0,\quad\forall\,x\in S\ \text{with}\ \Vert x\Vert=1.
\end{equation}
Namely, for $t=t_0$, the inequality \eqref{eq:claim2} follows from
the fact that the restriction of $B(t_0)$ to $S$ is of the form
$\langle(\mathrm{Id}+\overline K)\cdot,\cdot\rangle$ for some compact
operator $\overline K:S\mapsto S$. In this case, $c_0$ may be chosen to be the least
eigenvalue of $\mathrm{Id}+\overline K$. The continuity of $B$ concludes the
proof of the claim.

We set:
\begin{equation}\label{eq:4.2.4}
c_1=\inf_{\stackrel{y\in N}{\Vert y\Vert=1}}B'(t_0)[y,y]>0.
\end{equation}
Since $B$ is $C^1$, it is easy to see that, for $t$ sufficiently close 
to $t_0$, it is:
\begin{equation}\label{eq:4.2.5}
\phantom{\quad\forall\,y\in N,\ \Vert y\Vert=1.}
B(t)[y,y]\ge\frac12\,c_1\,(t-t_0),\quad\forall\,y\in N,\ \Vert y\Vert=1,
\end{equation}
so that $B(t)$ is positive 
definite on both $N$ and $S$ for $t$ sufficiently close to $t_0$. 
We want to show that, if $t>t_0$ is sufficiently close to $t_0$, then 
for all $x\in S\setminus\{0\}$
and $y\in N\setminus\{0\}$,
$B(t)$ is positive definite on the two dimensional subspace of $\mathcal H$ generated by
$x$ and $y$. By the positivity on $S$ and $N$, it suffices to prove that,
for $t>t_0$ is sufficiently close to $t_0$, the 
following inequality holds:
\begin{equation}\label{eq:schwartz}
B(t)[x,y]^2< B(t)[x,x]\cdot B(t)[y,y],
\end{equation}
for all $x\in S$, $y\in N$, $x,y\ne0$. Obviously, we can assume
$\Vert x\Vert=\Vert y\Vert=1$. As $B(t_0)$ vanishes on $N\times S$ and $B$ is of class
$C^1$, 
there exists $c_2>0$ such that, for all $t>t_0$ is sufficiently close to $t_0$, 
we have:
\begin{equation}\label{eq:second}
\big\vert B(t)[x,y]\big\vert\le c_2\cdot (t-t_0),
\end{equation}
for all $x\in S$, $y\in N$ with $\Vert x\Vert=\Vert y\Vert=1$.
By \eqref{eq:4.2.4}, \eqref{eq:4.2.5} and \eqref{eq:second},
for all $t>t_0$ is sufficiently close to $t_0$ we get:
\[B(t)[x,y]^2\le c_2^2\,(t-t_0)^2<\frac12 c_0\,c_1\,(t-t_0)\le B(t)[x,x]\cdot B(t)[y,y],\]
for all $x\in S$, $y\in N$ with $\Vert x\Vert=\Vert y\Vert=1$.
This yields \eqref{eq:schwartz} and concludes the first part of the proof.

For the general case, we use the spectral decomposition of $K$ to write 
an orthogonal decomposition $\mathcal H=S_+\oplus S_-\oplus N$, where
$B(t_0)$ is positive definite on $S_+$ and negative definite on $S_-$;
observe that $S_-$ is finite dimensional, and $n_-(B(t_0))=\mathrm{dim}(S_-)$.
Moreover, we write $N=N_+\oplus N_-$, where $B'(t_0)$ is positive definite
on $N_+$ and negative definite on $N_-$. We then apply the result proven 
in the first part of the proof to the restriction of $B(t)$ to $S_+\oplus N_+$ 
once, and again to the restriction of $-B(t)$ to $S_-\oplus
N_-$%
\footnote{observe that $S_-\oplus N_-$ has finite dimension, hence it is trivial that
the restriction of $-B(t)$ to $S_-\oplus N_-$ is represented by a compact perturbation
of a positive isomorphism, say the identity, and the first part of the proof applies.}.
 The conclusion follows by observing that $B(t)$ is positive definite
on $S_+\oplus N_+$ and negative definite on $S_-\oplus N_-$, which
implies that  $n_-(B(t))=
{\rm dim}(S_-\oplus N_-)$ for $t$ sufficiently close to $t_0$. Clearly, this also implies
that $B(t)$ is non degenerate.
\end{proof}
Although we will not need it, we observe that, for $t$ sufficiently close to $t_0$,
the bilinear map $\overline B(t)$ is actually strongly non degenerate, as it follows
easily from Fredholm's Alternative. We also observe that the assumption that
the bilinear map $\overline B(t_0)$ be represented by a compact perturbation 
of a positive operator cannot be removed from the statement of 
Proposition~\ref{thm:HSelementary}; it is easy to give examples
where the hypothesis is not satisfied and the thesis of Proposition~\ref{thm:HSelementary}
does not hold.
\begin{remdegenerate}\label{thm:remdegenerate}
It is important to emphasize that the conclusion of Proposition~\ref{thm:HSelementary}
does {\em not\/} hold if the assumption of nondegeneracy for the derivative $B'(t_0)$ is
not satisfied, and this is trivially checked. Besides, unless the Hilbert space
$\mathcal H$ is one-dimensional, it is very unlikely that the conclusion of
Proposition~\ref{thm:HSelementary} can be extended if one only makes
a non degeneracy assumption for some higher order derivative $B^{(k)}(t_0)$
on $\mathrm{Ker}(B(t_0))$; to understand this, we consider the following example.
%
Let $B_1(t)$ and $B_2(t)$ be the symmetric bilinear forms on $\R^2$ represented 
with respect to the canonical basis by  the following matrices:
\begin{equation}\label{eq:B1B2}
B_1(t)=\left(\begin{array}{cc}t^2&t\\ t& 1+t \end{array}\right),\quad
B_2(t)=\left(\begin{array}{cc}t^2&0\\ 0& 1  \end{array}\right).
\end{equation}
Clearly, $t_0=0$ is an isolated singularity for both $B_1$ ad $B_2$, and
$\mathrm{Ker}(B_1(0))=\mathrm{Ker}(B_2(0))=\R\cdot e_1$, where $e_1$ is the first
vector of the canonical basis of $\R^2$. The derivatives $B_1'(0)$ and $B_2'(0)$
vanish on $\R\cdot e_1$; moreover, the restrictions
of $B_1(t)$ and $B_2(t)$ on $\R\cdot e_1$ coincide for all $t$.
However, the change of value of the functions $n_-(B_1(t))$ and $n_-(B_2(t))$ 
passing from a negative to a positive value of $t$ is different:
\[\begin{split}
&n_-(B_1(t))=1,\quad\text n_-(B_2(t))=0,\quad\text{for}\ t<0,\\
&n_-(B_1(t))=0,\quad\text n_-(B_2(t))=0,\quad\text{for}\ t>0.
\end{split}\]
\end{remdegenerate}

\begin{remconstant}\label{thm:remconstant}
Observe that,
under the hypotheses of Proposition~\ref{thm:HSelementary}, 
if $B(t)$ is non degenerate for $t$ in some interval $I$, then
the function $i(t)=n_-(B(t))$ is constant on~$I$.
We also observe that Proposition~\ref{thm:HSelementary} can be 
applied to a {\em backwards reparameterization\/} of the curve $B(t)$ to
obtain information about the value of $n_-(B(t))$ for $t<t_0$ sufficiently close
to $t_0$. Namely, if one considers the curve of bilinear maps $S(t)=B(t_0-t)$,
we have $S(0)=B(t_0)$, $S'(0)=-B'(t_0)$, and the equality \eqref{eq:changen-} tells us
that, for $\tau>0$ sufficiently small, it is:
\begin{equation}\label{changen-back}
\begin{split}
n_-(\overline B(t_0-\tau))&=n_-(S(\tau))=n_-(\overline B(t_0))+n_-(-\widetilde B)=\\&
=n_-(\overline B(t_0))+n_+(\widetilde B).
\end{split}
\end{equation}
\end{remconstant}

We also have the following immediate corollary, which gives us a way to
compute the total change of index of a differentiable curve of
symmetric bilinear forms when passing through a degenerate instant:
\begin{corelementary}\label{thm:corelementary}
Let $B:[t_0-r,t_0+r]\mapsto\Bsym(\mathcal H,\R)$ and $\{\mathcal D_t\}_{t\in[t_0-r,t_0+r]}$
satisfy the same hypotheses of Proposition~\ref{thm:HSelementary}. Then, in the notations of
Proposition~\ref{thm:HSelementary},
for $\varepsilon>0$ small enough, we have:
\begin{equation}\label{eq:doublechangen-}
n_-(\overline B(t_0-\varepsilon))-n_-(\overline B(t_0+\varepsilon))=\sgn (\widetilde B).
\end{equation}
\end{corelementary}
\begin{proof}
Use Proposition~~\ref{thm:HSelementary} twice, once to $B\vert_{[t_0,t_0+r]}$
and once to a backwards  reparameterization of $B\vert_{[t_0-r,t_0]}$ 
(see Remark~\ref{thm:remconstant}).
\end{proof}
We conclude the section by showing a method that will be used later
to produce $C^1$-families of closed subspaces of a Hilbert space:
\begin{produce}\label{thm:produce}
Let $I\subset\R$ be an interval, $\mathcal H,\tilde{\mathcal H}$ be Hilbert
spaces and $F:I\mapsto\mathcal L(\mathcal H,\tilde{\mathcal H})$ be a $C^1$-map
such that each $F(t)$ is surjective. Then, the family $\mathcal D_t=\mathrm{Ker}(F(t))$
is a $C^1$-family of closed subspaces of $\mathcal H$.
\end{produce}
\begin{proof}
We exhibit local trivializations for the family $\{\mathcal D_t\}_{t\in I}$.
For $t=t_0\in I$, the map $F(t)$ maps the orthogonal complement
$\mathcal D_{t_0}^\perp$ isomorphically onto
$\tilde{\mathcal H}$; by continuity, this also holds for $t$ sufficiently close to $t_0$.
This implies that we have a direct sum decomposition $\mathcal H=\mathcal D_t\oplus
\mathcal D_{t_0}^\perp$ and the projection $\pi_t$ onto $\mathcal D_t$ is given by:
\[\pi_t=\mathrm{Id}-(F(t)\vert_{\mathcal D_{t_0}^\perp})^{-1}\circ F(t).\]
Obviously, $t\mapsto\pi_t$ is $C^1$. For $t$ sufficiently close to $t_0$, we 
define $\alpha(t)$ to be the inverse of the isomorphism:
\[(\pi_t\oplus\mathrm{Id}):\mathcal D_{t_0}\oplus\mathcal D_{t_0}^\perp\mapsto
\mathcal D_{t}\oplus\mathcal D_{t_0}^\perp.\]
Such a map $\alpha$ gives the required local trivialization for the
family $\{\mathcal D_t\}_{t\in I}$.
\end{proof}
\end{section}
\begin{section}{Morse--Sturm Systems and the Index Theorem for positive definite
metrics.}
\label{sec:applications}
Motivated by a geometric problem, we introduce a set of data $(g,R,P,S)$ for the Morse--Sturm 
problem as follows. Let's consider the system of differential equations in $\R^n$:
\begin{equation}\label{eq:MS}
\phantom{\quad t\in[0,1]}
J''(t)=R(t)[J(t)],\quad t\in[0,1]
\end{equation}
with initial conditions:
\begin{equation}\label{eq:IC}
J(0)\in P,\quad J'(0)+S[J(0)]\in P^\perp,
\end{equation}
where:
\begin{itemize}
\item $g$ is a (fixed) nondegenerate symmetric bilinear form on $\R^n$;
\smallskip

\item $R:[0,1]\mapsto\mathcal L(\R^n,\R^n)$ is a continuous map of {\em $g$-symmetric\/}
linear maps on $\R^n$, i.e., $g(R(t)[x],y)=g(x,R(t)[y])$ for all $x,y\in\R^n$;
\smallskip

\item $P$ is a subspace of $\R^n$ on which $g$ is non degenerate, and $P^\perp$ denotes
the orthogonal space of $P$ with respect to $g$;
\footnote{henceforth, the symbol $\perp$ will mean orthogonality with respect to $g$.}
\smallskip

\item $S:P\mapsto P$ is a $g$-symmetric linear map.
\end{itemize}
In some of the statements proven in this section, we will 
assume that $R$ is indeed a map of class $C^1$. Nevertheless, 
some perturbation arguments presented in the next section will 
allow us to prove our main results in the general case of a continuous
map $R$.

A solution for the differential equation \eqref{eq:MS} satisfying the initial
conditions \eqref{eq:IC} will be called a $(P,S)$-solution; we denote by
$\mathbb J$ the set of all $(P,S)$-solutions:
\begin{equation}\label{eq:defbbJ}
\mathbb J=\Big\{J:[0,1]\mapsto\R^n:J\ \text{satisfies \eqref{eq:MS} 
and \eqref{eq:IC}}\Big\}.
\end{equation}
Observe that $\mathbb J$ is an $n$-dimensional vector space.
For all $t\in [0,1]$, we define $\mathbb J[t]$ by:
\begin{equation}\label{eq:defbbJt}
\mathbb J[t]=\Big\{J(t):J\in\mathbb J\Big\},
\end{equation}
and we say that $t_0\in \,]0,1]$ is a {\em $(P,S)$-focal instant\/} if
there exists a non zero $J\in\mathbb J$ such that $J(t_0)=0$.
Clearly, this is equivalent to requiring that \/ $\mathbb J[t_0]\ne \R^n$.
The {\em multiplicity\/} $\mu(t_0)$ of a $(P,S)$-focal instant $t_0$ is the codimension
of $\mathbb J[t_0]$ in $\R^n$, or equivalently, the dimension of $\mathbb J[t_0]^\perp$.
The {\em signature\/} $\sgn(t_0)$ of $t_0$ is defined as the signature of
the restriction of the bilinear form $g$ to the space $\mathbb J[t_0]^\perp$:
\begin{equation}\label{eq:defsignature}
\sgn(t_0)=\sgn\left(g\big\vert_{\mathbb J[t_0]^\perp}\right).
\end{equation}
The $(P,S)$-focal instants coincide with the set of zeroes of the function
$r(t)={\mathrm det}(J_1(t),J_2(t),\ldots,J_n(t))$, where $J_1,\ldots,J_n$
is a basis of $\mathbb J$. If $R(t)$ is real analytic, then also $r(t)$ is
real analytic on $[0,1]$, hence its zeroes are isolated (observe that
$r(t)$ cannot be identically zero, see Proposition~\ref{thm:discrete}). 
In \cite[Proposition~2.5.1]{MPT} some sufficient conditions
for the discreteness of the $(P,S)$-focal instants are given. More precisely,
the following result is proven:
\begin{discrete}\label{thm:discrete}
Let $t_0$ be a $(P,S)$-focal instant. If $g$ is non degenerate on $\mathbb J[t_0]$
(or equivalently on $\mathbb J[t_0]^\perp$) then there are no other $(P,S)$-focal
instants in some neighborhood of $t_0$. Moreover, there are no $(P,S)$-focal instants
in some neighborhood of $t_0=0$.
\qed
\end{discrete}
A proof of Proposition~\ref{thm:discrete} can also be deduced
from some results that will be presented in the rest of this 
section (see Remark~\ref{thm:remdiscrete}).\smallskip

An easy calculation shows that, for $J_1,J_2\in\mathbb J$, the following equality
holds:
\begin{equation}\label{eq:symmetry}
\phantom{,\quad\forall\,t\in[0,1].}
g(J_1'(t),J_2(t))=g(J_1(t),J_2'(t)),\quad\forall\,t\in[0,1].
\end{equation}
Namely, we use \eqref{eq:MS} to show that the difference $g(J_1',J_2)-g(J_1,J_2')$ is constant,
and  \eqref{eq:IC} to see that this constant is zero. Formula \eqref{eq:symmetry}
and an easy dimension counting argument shows that, for $t\in[0,1]$:
\begin{equation}\label{eq:Jperp'}
\mathbb J[t]^\perp=\Big\{J'(t):J\in\mathbb J,\ J(t)=0\Big\}.
\end{equation}
Namely, from \eqref{eq:symmetry} it follows easily the inclusion of
the term on the right hand side into $\mathbb J[t]^\perp$; conversely, 
it is easy to see that the dimension of the space on the right hand
side of \eqref{eq:Jperp'} is equal to $\mu(t)=\mathrm{dim}(\mathbb J[t]^\perp)$,
which proves \eqref{eq:Jperp'}.

Moreover, we introduce the following analytical framework.

Let $H^1([a,b],\R^m)$ denote the Sobolev space of all absolutely continuous
$\R^m$-valued maps on $[a,b]$ with square integrable derivative;
$H^1_P([a,b],\R^m)$ will denote the subspace of $H^1([a,b],\R^m)$ 
consisting of those $V$ 
such that $V(a)\in P$ and $V(b)=0$. Moreover, $H^1_0([a,b],\R^m)$ is the subspace
of $H^1([a,b],\R^m)$ given by the $V$'s such that $V(a)=V(b)=0$.

For $t\in\,]0,1]$, we set $\mathcal H_t=H^1_P([0,t],\R^n)$ and $\mathcal H=\mathcal H_1$;
we define the isomorphisms 
\begin{equation}\label{eq:defvarphit}
\varphi_t:\mathcal H\mapsto\mathcal H_t, \quad\text{with}\quad 
\varphi_t(\hat V)(s)=V(s)=\hat V\left(\frac st\right),\quad s\in [0,t].\end{equation} 
For each $t\in\,]0,1]$,  we introduce the {\em index form\/} $I_t$ on $\mathcal H_t$, 
which is the symmetric bilinear form given by:
\begin{equation}\label{eq:defIt}
\begin{split}
I_t(V&,W)=\\
=&\int_0^t \left[g(V'(s),W'(s))+g(R(s)[V(s)],W(s))\right]\;{\mathrm d}s- g(S[V(0)],W(0)).
\end{split}
\end{equation}
\begin{remId+K}\label{thm:remId+K}
If $g$ is positive definite, then one can consider the following Hilbert space
inner product on $\mathcal H_t$:
\[\langle V,W\rangle_{\mathcal H_t}=\int_0^tg(V'(s),W'(s))\,\mathrm ds.\]
The bilinear form $I_t$ is written as the sum of $\langle\cdot,\cdot\rangle_{\mathcal H_t}$
and a bilinear form which is continuous with respect to the $C^0$-topology. 
By the compact embedding of $H^1([0,t],\R^n)$ in $C^0([0,t],\R^n)$ (see~\cite{Brezis}),
one obtains immediately that $I_t$ is of the form $\langle
(\mathrm{Id}+K)\,\cdot\,,\,\cdot\,\rangle_{\mathcal H_t}$ for some compact operator $K$ on $\mathcal H_t$.
\end{remId+K}
Finally, for all $t\in\,]0,1]$, let $\hat I_t$ be the symmetric bilinear form on $\mathcal H$ obtained
by the pull-back of $I_t$ by $\varphi_t$, namely:
\begin{equation}
\hat I_t=I_t(\varphi_t\cdot,\varphi_t\cdot).
\end{equation}
Explicitly, for $\hat V,\hat W\in\mathcal H$ we have:
\begin{equation}\label{eq:explicitBt}\begin{split}
\hat I_t(\hat V&,\hat W)=\\
=&\int_0^t \left[\frac1{t^2}\,g\left(\hat V'\left(\frac st\right),\hat W'\left(\frac
st\right)\right)+g\left(R(s)\left[\hat V\left(\frac st\right)\right],\hat W\left(\frac
st\right)\right)\right]\;{\mathrm d}s\\ &- g(S[\hat V(0)],\hat W(0)).
\end{split}
\end{equation}
Integration by parts in \eqref{eq:defIt} and the Fundamental Lemma of Calculus of Variations
show that 
\begin{equation}\label{eq:KerIt}
\mathrm{Ker}(I_t)=\Big\{J\big\vert_{[0,t]}:J\in\mathbb J,\ J(t)=0\Big\};
\end{equation}
from \eqref{eq:Jperp'} and \eqref{eq:KerIt} for each $t\in\,]0,1]$ we then get an isomorphism
\begin{equation}\label{eq:defpsit}
\begin{split}
\psi_t:\mathrm{Ker}(I_t)&\longmapsto\;\; \mathbb J[t]^\perp\\
V\quad&\longmapsto V'(t).
\end{split}
\end{equation}
We set
\begin{equation}\label{eq:defNt}
\mathcal N_t=\mathrm{Ker}(\hat I_t)\subset\mathcal H;
\end{equation}
obviously, $\varphi_t$ gives an isomorphism between $\mathrm{Ker}(I_t)$ and $\mathcal N_t$.
\begin{BC1}\label{thm:BC1}
Suppose that $R$ is a map of class $C^1$. Then,
the map \[]0,1]\ni t\mapsto \hat I_t\in\mathrm B_{\mathrm{sym}}(\mathcal H,\R)\] is of class $C^1$.
Moreover, the map $]0,1]\ni t\mapsto C_t=t\cdot \hat I_t$ has a $C^1$-extension to $[0,1]$, with
\begin{equation}\label{eq:defC0}
C_0(\hat V,\hat W)=\int_0^1g(\hat V'(u),\hat W'(u))\;\mathrm du,\quad \hat V,\hat W\in\mathcal H.
\end{equation}
\end{BC1}
\begin{proof}
Substituting $u=\frac st$ in \eqref{eq:explicitBt}, we get the following expression
for $\hat I_t$:
\begin{equation}\label{eq:Bt2}
\begin{split}
\hat I_t(\hat V,\hat W)=&\int_0^1\left[\frac1t\, g(\hat V'(u),\hat W'(u))+t g(R(tu)[\hat V(u)],\hat W(u))
\right]\;\mathrm du\\&-g(S[\hat V(0)],\hat W(0)).
\end{split}
\end{equation}
Differentiating \eqref{eq:Bt2} with respect to $t$ we get:
\begin{equation}\label{eq:Bt'}
\begin{split}
\frac{\mathrm d}{\mathrm dt}\,\hat I_t(\hat V,\hat W)=&\int_0^1\left[-\frac1{t^2}\,
g(\hat V'(u),\hat W'(u))+g(R(tu)[\hat V(u)],\hat W(u))\right]\;\mathrm du
\\&+t\int_0^1 u\,g(R'(tu)[\hat V(u)],\hat W(u))\;\mathrm du.
\end{split}
\end{equation}
We now apply Lemma~\ref{thm:C1} to $F(t)=\hat I_t$, $G(t)$ is the right hand side of
equality \eqref{eq:Bt'}, $E=\mathrm B_{\mathrm{sym}}(\mathcal H,\R)$, $E_0=\R$
and $\Phi=\{\phi_{\hat V,\hat W}: \hat V,\hat W\in\mathcal H\}$, where
\[ 
\phi_{\hat V,\hat W}(B)=B(\hat V,\hat W),\quad B\in\mathrm B_{\mathrm{sym}}(\mathcal H,\R).\]
It is easy to check that $G$ is continuous, by the continuity of $R$ and $R'$, and clearly
$\Phi$ is separating for $E$, which concludes the first part of the proof.

From \eqref{eq:Bt2} we compute easily:
\begin{equation}\label{eq:tBt}
\begin{split}
C_t(\hat V,\hat W)=&\int_0^1\left[g(\hat V'(u),\hat W'(u))+t^2 g(R(tu)[\hat
V(u)],\hat W(u))
\right]\;\mathrm du\\&-t\cdot g(S[\hat V(0)],\hat W(0)),
\end{split}
\end{equation}
for all $t\in[0,1]$. Its regularity is established analogously
applying Lemma~\ref{thm:C1}.
\end{proof}
We have the following immediate Corollary:
\begin{corpositive}\label{thm:corpositive}
For $t>0$ small enough, $I_t$ is strongly non degenerate on $\mathcal H_t$.
Moreover, if $g$ is positive definite in $\R^n$, then $I_t$ is
positive definite for $t$ small enough.
\end{corpositive}
\begin{proof}
For $t>0$, $I_t$ is strongly non degenerate (positive) 
if and only if $C_t$ is strongly
non degenerate (positive). From \eqref{eq:defC0}, $C_0$ is strongly non degenerate
because $g$ is non degenerate; by continuity, $C_t$ is also strongly non degenerate
for $t>0$ small enough.

If $g$ is positive definite, then $C_0$ is a Hilbert space inner product, 
and therefore it is positive definite and away
from $0$. By continuity, $C_t$ is positive definite for $t$ small enough.
\end{proof}
We now pass to the study of the signature of $\hat I'(t)$ on $\mathcal N_t$.
For this, we consider the push-forward of $\hat I'(t)$ through the isomorphism:
\[\psi_t\circ\varphi_t:\mathcal N_t\longmapsto\mathbb J[t]^\perp \]
given by the composition:
\[\hat I'(t)\left((\psi_t\circ\varphi_t)^{-1}\cdot,(\psi_t\circ\varphi_t)^{-1}\cdot\right),\]
where the maps $\varphi_t$ and $\psi_t$ are defined in \eqref{eq:defvarphit} and 
\eqref{eq:defpsit}.

We have the following:
\begin{central}\label{thm:central}
Suppose that $R$ is a map of class $C^1$.
For $t\in\,]0,1]$, the isomorphism $\psi_t\circ\varphi_t$ carries the restriction
of $\hat I'(t)$ to $\mathcal N_t$ into the restriction of $-g$ to $\mathbb J[t]^\perp$.
\end{central}
\begin{proof}
Let $t\in \,]0,1]$ and $\hat V,\hat W\in\mathcal N_t$ be fixed; 
observe that $\hat V$ and $\hat W$
are maps of class $C^3$, because they are affine
reparameterizations of solutions to \eqref{eq:MS};
they satisfy the following differential equations:
\begin{equation}\label{eq:repED}
\frac1{t^2}\,\hat V''\left(\frac st\right)=R(s)\,\left[\hat V\left(\frac st\right)\right],\quad
\frac1{t^2}\,\hat W''\left(\frac st\right)=R(s)\,\left[\hat W\left(\frac st\right)\right],
\qquad s\in [0,t]. 
\end{equation}
We differentiate \eqref{eq:explicitBt} with respect to $t$ 
and, observing that $\hat V(1)=\hat W(1)=0$, we obtain:
\begin{equation}\label{eq:infernal}
\begin{split}
\frac{\mathrm d}{\mathrm dt}\,\hat I_t(\hat V,\hat W)=&\frac1{t^2}\,g(V'(1),W'(1))-
\int_0^t\frac2{t^3}\,g\left(\hat V'\left(\frac st\right),\hat W'\left(\frac st\right)\right)
\;\mathrm ds\\ -&\int_0^t\frac s{t^4}\left[g\left(\hat V''\left(\frac st\right),\hat W'\left(\frac
st\right)\right)+g\left(\hat V'\left(\frac st\right),\hat W''\left(\frac
st\right)\right)\right]\;\mathrm ds\\-
\int_0^t\frac s{t^2}&\left[g\left(R(s)\hat  V'\left(\frac st\right),\hat W\left(\frac
st\right)\right)+g\left(R(s)\hat V\left(\frac st\right),\hat W'\left(\frac
st\right)\right)\right]
\;\mathrm ds.
\end{split}
\end{equation}
Using \eqref{eq:repED}, we eliminate from \eqref{eq:infernal} the terms involving
the operator $R$, and we get:
\begin{equation}\label{eq:afterinfernal}
\begin{split}
\frac{\mathrm d}{\mathrm dt}\,\hat I_t(\hat V,\hat W)=&\frac1{t^2}\,g(V'(1),W'(1))-2
\int_0^t\frac{\mathrm d}{\mathrm ds}\left[\frac s{t^3}g\left(\hat V'\left(\frac
st\right),\hat W'\left(\frac st\right)\right)\right]\;\mathrm ds=
\\&=-\frac 1{t^2}\,g(\hat V(1),\hat W(1))=-g(\psi_t\circ\varphi_t(\hat V),\psi_t\circ\varphi_t(\hat W)).
\end{split}
\end{equation}
This concludes the proof.
\end{proof}
\begin{remdiscrete}\label{thm:remdiscrete}
If $t_0$ is a $(P,S)$-focal instant for which $g$ is non degenerate on $\mathbb J[t_0]$,
then Proposition~\ref{thm:HSelementary} and  Proposition~\ref{thm:central} imply that
$\hat I_t$, and hence $I_t$, is non degenerate for $t\ne t_0$ sufficiently close to $t_0$.
Moreover, by Corollary~\ref{thm:corpositive} there are no $(P,S)$-focal instants near
$t=0$. So, we obtain an alternative proof of Proposition~\ref{thm:discrete}.
\end{remdiscrete}
As a corollary to Proposition~\ref{thm:HSelementary} and Proposition~\ref{thm:central},
we obtain the classical Morse--Sturm Oscillation Theorem:
\begin{RiemannianMorse}\label{thm:RiemannianMorse}
Suppose that $R$ is a map of class $C^1$.
If $g$ is positive definite in $\R^n$, then the following equality holds:
\begin{equation}\label{eq:MSO}
n_-(I_1)=\sum_{t\in\,]0,1[} \mu(t).
\end{equation}
\end{RiemannianMorse}
\begin{proof}
Let $t_0\in\,]0,1]$ be fixed. By Remark~\ref{thm:remId+K}, $\hat I_{t_0}$ is represented
by a compact perturbation of the identity map with respect to some suitably chosen
Hilbert space inner product on $\mathcal H$. By Proposition~\ref{thm:BC1}, $\hat I$ is
of class $C^1$, and we are under the hypotheses of Proposition~\ref{thm:HSelementary}.
If $t_0<1$, applying Corollary~\ref{thm:corelementary} and Proposition~\ref{thm:central}, 
we obtain that the integer valued function $i(t)= n_-(\hat I_t)$ is constant around $t_0$ if $t_0$ is
not a $(P,S)$-focal instant, whereas it has a {\em jump\/} of exactly $\mu(t_0)$
at $t_0$ if $t_0$ is a $(P,S)$-focal. If $t$ is small enough, 
by Corollary~\ref{thm:corpositive}, it is $n_-(\hat I_t)=0$, and this concludes
the proof in the case that $t_0=1$ is not $(P,S)$-focal.

\noindent\ 
Applying Proposition~\ref{thm:HSelementary} to backwards reparameterizations
of $\hat I_t$ (see Remark~\ref{thm:remconstant}), we see that the map $i(t)$ is indeed a {\em
left-continuous\/} function on $]0,1]$, and therefore $n_-(\hat I_1)=n_-(\hat I_{1-\varepsilon})$ for
$\varepsilon>0$ small enough. With this observation the proof is concluded.
\end{proof}
\end{section}
\begin{section}{The Index Theorem for non positive definite metrics}
\label{sec:nonpositive}
In this section we aim at a generalization of the result of Corollary~\ref{thm:RiemannianMorse}
to the case of non positive definite metrics $g$. As we have observed, for a general
metric $g$ the left-hand side of the equality \eqref{eq:MSO} is infinite; on the other
hand, the sum appearing in the right-hand side of \eqref{eq:MSO} may lose sense, due the 
fact that there may be an infinity of  focal instants.

For the beginning, we will consider only the case of Morse--Sturm systems
having a finite number of $(P,S)$-focal instants. We will see
that this assumption holds {\em generically}, i.e., for almost all choices of
the data $R,P,S$ in \eqref{eq:MS} and \eqref{eq:IC}. The conclusion for the
general case will be obtained by perturbation arguments, discussed in Section~\ref{sec:nondeg}.
As to the finiteness of the index,
we want to consider the restriction of $I_1$ to a suitable subspace
$\mathcal K$ of $\mathcal H$ that ought to be {\em small enough\/} to yield finiteness
of the index, 
but {\em large enough\/} to retain the relevant information about the
differential problem. Actually, in order to use the techniques 
of Section~\ref{sec:applications} to
compute the evolution of the index function $i(t)$, we need to 
determine a whole family $\mathcal K_t$
of subspaces of $\mathcal H_t$ with the required properties.

Having a concrete example in mind, we axiomatize the
following set of properties for the family $\mathcal K_t$.
\begin{defadmissible}\label{thm:defadmissible}
For each $t\in\,]0,1]$, let $\mathcal K_t$ be a closed subspace
of $\mathcal H_t$ and let $\hat{\mathcal K}_t=\varphi_t^{-1}(\mathcal K_t)$.
The family $\{\mathcal K_t\}_{t\in\,]0,1]}$ is called an {\em admissible
family of subspaces\/} for the Morse--Sturm Problem \eqref{eq:MS} and
\eqref{eq:IC} if the following conditions are satisfied:
\begin{enumerate}
\item\label{itm:C1} the family $\{\hat{\mathcal K}_t\}$ admits an extension to
$t=0$, denoted by $\hat{\mathcal K}_0$, which makes it a $C^1$-family
of closed subspaces on the interval $[0,1]$;
\item\label{itm:poscomp} for $t\in\,]0,1]$, the restriction of the index form
$I_t$ to $\mathcal K_t$ is represented by a linear operator which is the 
sum of a positive self-adjoint isomorphism of $\mathcal K_t$ and a compact
(self-adjoint) operator on $\mathcal K_t$;
\item\label{itm:poscomp0}  the restriction of the bilinear form $C_0$ (see formula~\eqref{eq:defC0})
to $\hat{\mathcal K}_0$ is non degenerate, and it is
represented by the sum of a positive self-adjoint 
isomorphism and a compact (self-adjoint) operator on $\hat{\mathcal K}_0$;
\item\label{itm:kernel} for $t\in\,]0,1]$, the kernel of the restriction of $I_t$ to
$\mathcal K_t$ is equal to the kernel of $I_t$ in $\mathcal H_t$ (see formula~\eqref{eq:KerIt}).
\end{enumerate}
\end{defadmissible}
The condition~\ref{itm:poscomp} of Definition~\ref{thm:defadmissible} implies that,
for each $t\in\,]0,1]$ there exists a Hilbert space inner product on $\mathcal K_t$
under which the bilinear form $I_t$ is represented by a compact perturbation
of the identity map on $\mathcal K_t$. By condition~\ref{itm:poscomp0}, the same 
is true for the bilinear map $C_0$ on $\hat{\mathcal K}_0$. In particular, by the 
condition~\ref{itm:C1} and by Proposition~\ref{thm:BC1}, we are
allowed to use the result of Proposition~\ref{thm:HSelementary} and of 
Corollary~\ref{thm:corelementary} to the bilinear forms $\hat I_t$ and $C_t$ on 
$\hat{\mathcal K}_t$. Observe that the hypothesis~\ref{itm:hp2.5.3} of Proposition~\ref{thm:HSelementary}
for the family of closed subspaces $\mathcal K_t$
is satisfied thanks to the axioms~\ref{itm:poscomp0} and \ref{itm:kernel} 
of Definition~\ref{thm:defadmissible}.

The axioms satisfied by an admissible family of subspaces for the
Morse--Sturm problem constitute the hypotheses of a generalization
of Corollary~\ref{thm:RiemannianMorse}. Recalling the definition 
\eqref{eq:defsignature} of the signature $\sgn(t)$ of a $(P,S)$-focal instant $t$,
we prove the following:
\begin{indexthm}\label{thm:indexthm}
Let $\{\mathcal K_t\}_{t\in\,]0,1]}$ be an admissible
family of subspaces for the Morse--Sturm problem \eqref{eq:MS} and
\eqref{eq:IC}, with $R$ of class $C^1$, and assume that the restriction
of $g$ to $\mathbb J[t]$ is non degenerate for all $t\in\,]0,1]$. Then, we have
the following equality:
\begin{equation}\label{eq:indextheorem}
n_-(I_1\vert_{\mathcal K_1})=n_-(C_0\vert_{\hat{\mathcal
K}_0})+\sum_{t\in\,]0,1[}\sgn(t)-n_-(g\vert_{\mathbb J[1]^\perp}).
\end{equation}
\end{indexthm}
\begin{proof}
As in the proof of Corollary~\ref{thm:RiemannianMorse},
we study the {\em evolution\/} of the function
$i(t)=n_-(\hat I_t\vert_{\hat{\mathcal K}_t})$ when
$t$ runs from $0$ to $1$; observe that $i(1)=n_-(\hat I_1\vert_{\hat{\mathcal
K}_1})=n_-(I_1\vert_{\mathcal K_1})$. Observe that, by the axiom~\ref{itm:kernel}
of Definition~\ref{thm:defadmissible}, the $(P,S)$-focal instants coincide
with the instants $t$ where $\hat I_t$ is degenerate on $\hat{\mathcal K}_t$.

By Proposition~\ref{thm:discrete} (see also Remark~\ref{thm:remdiscrete}), there is only a
finite number of $(P,S)$-focal instants, hence, by Remark~\ref{thm:remconstant}, 
$i$ is piecewise constant on $]0,1]$. Namely, $i$ is constant on any interval that
does not contain $(P,S)$-focal instants.

Since  $n_-(C_t\vert_{\hat{\mathcal K}_t})=n_-(\hat I_t\vert_{\hat{\mathcal K}_t})$ for $t>0$,
by the non degeneracy of $C_0$ on $\hat{\mathcal K}_0$ and Remark~\ref{thm:remconstant},
$i(t)=n_-(C_0\vert_{\hat{\mathcal K}_0})$ for $t>0$ sufficiently small. 

When $t$ {\em passes\/} through a $(P,S)$-focal instant $t_0\in\,]0,1[$, by
Corollary~\ref{thm:corelementary} and by Proposition~\ref{thm:central} the
jump of the function $i$ is equal to the signature  $\sgn(t_0)$.

\noindent\ \ 
Finally, applying Proposition~\ref{thm:HSelementary} to a backwards reparameterization
of $\hat I_t$ around $t_0=1$ (see Remark~\ref{thm:remconstant}), by Proposition~\ref{thm:central} for
$t<1$ sufficiently close to $1$ we have $i(t)-i(1)=n_-(g\vert_{\mathbb J[1]^\perp})$, which concludes
the proof.
\end{proof}
We have observed in the proof of Corollary~\ref{thm:RiemannianMorse}
that the index function $i(t)$ is {\em left-continuous\/}
under the positivity assumption for $g$. We emphasize that, as it
was clear in the above proof, this property fails when $g$ is non positive. 
As a consequence of this lack of continuity, when comparing with the
Riemannian Index Theorem, in the right hand side
of equality \eqref{eq:indextheorem}  we get the extra term 
$n_-(g\vert_{\mathbb J[1]^\perp})$ which is non zero when $t_0=1$ is $(P,S)$-focal.

Another remarkable phenomenon that appears in the case of non positive
definite metrics is the presence of the term $n_-(C_0\vert_{\hat{\mathcal
K}_0}) $ in the equality \eqref{eq:indextheorem}, which is the {\em initial value\/}
of the index function $i(t)$. As we saw in the proof of Corollary~\ref{thm:RiemannianMorse}, 
for positive definite metrics, such initial value is zero.  
\smallskip

We now present a concrete example of the above situation. We will assume 
throughout the rest of this section that $n_-(g)=1$
and that the differential equation~\eqref{eq:MS} admits a solution
$Y:[0,1]\mapsto\R^n$ with the property that $g(Y,Y)<0$ on $[0,1]$:
\begin{equation}\label{eq:defY}
Y''=RY,\quad\text{and}\quad g(Y,Y)<0.
\end{equation}

We fix one such solution $Y$ and we consider the following one-parameter family
of positive definite inner products in $\R^n$:
\begin{equation}\label{eq:defgr}
\gr_t(v,w)=g(v,w)-2\,\frac{g(v,Y(t))g(w,Y(t))}{g(Y(t),Y(t))},
\quad\forall\,t\in[0,1],\ v,w\in\R^n.
\end{equation}
Observe that, for all $t\in[0,1]$, $\gr_t(v,w)$ coincides with $g(v,w)$ if either
$v$ or $w$ is orthogonal to $Y(t)$, and $\gr_t(Y(t),Y(t))=-g(Y(t),Y(t))$. The
formula that gives $g$ in terms of $\gr_t$ is similar:
\begin{equation}\label{eq:ggr}
g(v,w)=\gr_t(v,w)-2\,\frac{\gr_t(v,Y(t))\gr_t(w,Y(t))}{\gr_t(Y(t),Y(t))}.
\end{equation}
For all $t\in\,]0,1]$, we consider the following subspace of $\mathcal H_t$:
\begin{equation}\label{eq:defcalKt}
\mathcal K_t=\Big\{V\in\mathcal H_t:g(V',Y)-g(V,Y')\equiv C_V\ \text{(constant)}\Big\}.
\end{equation}
We claim that $\mathcal K_t$ is an {\em admissible
family of subspaces\/} for the Morse--Sturm Problem \eqref{eq:MS} and
\eqref{eq:IC}, and we take the rest of this section to prove the claim.

As in Definition~\ref{thm:defadmissible}, for $t\in\,]0,1]$ we set $\hat{\mathcal
K}_t=\varphi_t^{-1}(\mathcal K_t)$; explicitly, we have:
\begin{equation}\label{eq:defhatcalKt}
\hat{\mathcal K}_t=\Big\{\hat V\in\mathcal H: g(\hat V'(u),\hat
Y_t(u))-  g(\hat V(u),\hat Y_t'(u)))\equiv \text{const.}\Big\},
\end{equation}
where $\hat Y_t(u)=Y(t\cdot u)$ for $u\in[0,1]$. We observe that formula
\eqref{eq:defhatcalKt} makes sense also for $t=0$, where $\hat Y_0$ is the
constant vector $Y(0)$:
\begin{equation}\label{eq:defcalK0}
\hat{\mathcal K}_0=\Big\{\hat V\in\mathcal H:g(\hat V'(u),Y(0))\equiv\text{const.}\Big\}.
\end{equation}

Let $\widetilde{\mathcal H}$ denote the Hilbert space given  by the quotient
$L^2([0,1],\R)/\mathfrak C$, where $\mathfrak C$ denotes the subspace of constant
functions. For $t\in[0,1]$, $\hat{\mathcal K}_t$ is the kernel of the bounded
linear map $F_t:\mathcal H\mapsto\widetilde{\mathcal H}$ given by:
\begin{equation}\label{eq:dadiff}
\begin{split}
F_t(\hat V)(u)&=g(\hat V'(u),\hat
Y_t(u))-  g(\hat V(u),\hat Y_t'(u))+\mathfrak C=\\
&=g(\hat V'(u),Y(tu))- t\cdot g(\hat V(u),Y'(tu))+\mathfrak C.
\end{split}
\end{equation}
\begin{FtC1}\label{thm:FtC1} The map
$[0,1]\ni t\mapsto F_t\in \mathcal L(\mathcal H,\widetilde{\mathcal H})$ is of class $C^1$.
\end{FtC1}
\begin{proof}
We formally differentiate \eqref{eq:dadiff}, obtaining:
\begin{equation}\label{eq:formal}
\begin{split}
F'_t(V)(u)=\;&u\, g(\hat V'(u),Y'(tu))-g(\hat V(u),Y'(tu))+\\&-tu\,g(\hat V(u),Y''(tu))+
\mathfrak C.
\end{split}
\end{equation}
Using the fact that $Y$ is of class $C^2$, it is easily seen that formula \eqref{eq:formal}
defines a continuous curve in $\mathcal L(\mathcal H,\widetilde{\mathcal H})$.
We now use Lemma~\ref{thm:C1} by considering $\Phi$ to be the set of
{\em evaluations\/} at fixed vectors $\hat V\in\mathcal H$; the conclusion
will follow once we prove that the map $t\mapsto F_t(\hat V)\in\widetilde{\mathcal H}$
is of class $C^1$ for all $\hat V\in\mathcal H$, and that  its derivative is given by 
\eqref{eq:formal}.

Let $C^1([0,1],\R^n)$ be the Banach space of $\R^n$-valued $C^1$-maps on $[0,1]$;
we define the following bounded linear operator $\sigma:C^1([0,1],\R^n)\mapsto\widetilde{\mathcal H}$
by:
\begin{equation}\label{eq:defsigma}
\sigma(\mathcal Y)(u)=g(\hat V'(u),\mathcal Y(u))-g(\hat V(u),\mathcal Y'(u)).
\end{equation}
We observe that the map $t\mapsto F_t(\hat V)$ is given by the composition of $\sigma$ and 
the map 
\begin{equation}\label{eq:mapY}
t\mapsto \hat Y_t\in C^1([0,1],\R^n).
\end{equation}
It remains to show that the map \eqref{eq:mapY} is of class $C^1$. This is again an easy consequence
of Lemma~\ref{thm:C1}, where $\Phi$ is the set of evaluations at fixed instants $u\in[0,1]$.
\end{proof}
The next step towards our goal is to prove the surjectivity of $F_t$. We introduce the
subspaces $\mathcal S_t\subset\mathcal H_t$ and $\hat{\mathcal S}_t\subset\mathcal H$:
\begin{equation}\label{eq:defSt}
\begin{split}
&\mathcal S_t=\Big\{f\cdot Y\vert_{[0,t]}:f\in H^1_0([0,t],\R)\Big\},\quad t\in\,]0,1],\\
&\hat{\mathcal S}_t=\Big\{\hat f\cdot \hat Y_t:\hat f\in H^1_0([0,1],\R)\Big\},
\quad t\in[0,1].
\end{split}
\end{equation}
Observe that, for $t\in\,]0,1]$, $\mathcal S_t=\varphi_t(\hat{\mathcal S}_t)$. 
We show now that $F_t(\hat{\mathcal S}_t)=\widetilde{\mathcal H}$:
\begin{surjective}\label{thm:surjective}
For all $t\in[0,1]$, the restriction of $F_t$ to $\hat{\mathcal S}_t$ is
surjective.
\end{surjective}
\begin{proof} For $\hat f\in H^1_0([0,1],\R)$, we compute :
\[
F_t(\hat f\cdot \hat Y_t)=\hat f'\cdot g(\hat Y_t,\hat Y_t)+\mathfrak C.
\]
Hence, for the proof we need to show that, given $h\in L^2([0,1],\R)$ there
exists $c\in\R$ and $\hat f\in H^1_0([0,1],\R)$ such that the following
differential equation is satisfied:
\[\hat f'=\frac{h+c}{g(\hat Y_t,\hat Y_t)}.\]
It suffices to take:
\[c=-\left(\int_0^1\frac{\mathrm dr}{g(\hat Y_t,\hat Y_t)}\right)^{-1}\int_0^1
\frac h{g(\hat Y_t,\hat Y_t)}\;\mathrm dr, \quad\text{and}\quad \hat f(u)=\int_0^u\frac {h+c}{g(\hat
Y_t,\hat Y_t)}\;\mathrm dr.\]
Observe that the above formulas make  sense because $g(\hat Y_t,\hat Y_t)<0$.
\end{proof}
\begin{corC1}\label{thm:corC1}
$\{\hat{\mathcal K}_t\}_{t\in[0,1]}$ is a $C^1$-family of closed subspaces of
$\mathcal H$.
\end{corC1}
\begin{proof}
It follows directly from Lemma~\ref{thm:produce}, Lemma~\ref{thm:FtC1} 
and Lemma~\ref{thm:surjective}.
\end{proof}
\begin{corsoma}\label{thm:corsoma}
For $t\in\,]0,1]$, $\mathcal H_t=\mathcal K_t+\mathcal S_t$; moreover,
$\mathcal H=\hat{\mathcal K}_0+\hat{\mathcal S}_0$.
\end{corsoma}
\begin{proof}
By Lemma~\ref{thm:surjective}, 
an easy linear algebra argument shows that, for $t\in[0,1]$, $\mathcal
H=\hat{\mathcal K}_t +\hat{\mathcal S}_t$. For $t\in\,]0,1]$ we apply the isomorphism $\varphi_t$ and
we get the conclusion.
\end{proof}
Although we will not need it, we emphasize that the sums in the statement
of Corollary~\ref{thm:corsoma} are direct. As a matter of facts, we now
prove that the above sums are orthogonal with respect to the bilinear forms
$I_t$ and $C_0$, respectively.
\begin{orthogonal}\label{thm:orthogonal}
For all $t\in\,]0,1]$, the spaces ${\mathcal K}_t$ and ${\mathcal S}_t$ are
orthogonal with respect to the bilinear form $I_t$;
moreover, the spaces $\hat{\mathcal K}_0$ and $\hat{\mathcal S}_0$ are orthogonal
with respect to $C_0$.
\end{orthogonal}
\begin{proof}
Let $V\in\mathcal K_t$ and $f\cdot Y\in\mathcal S_t$ be fixed, with $f(0)=f(t)=0$. 
From \eqref{eq:defIt}, \eqref{eq:defY} and \eqref{eq:defcalKt}, we compute using integration
by parts as follows:
\begin{equation}
\begin{split}
I_t(V,fY)&=\int_0^t\left[f'g(V',Y)+f\,g(V',Y')+f\,g(RV,Y)\right]\;\mathrm ds=\\
&=\int_0^t\left[f'C_V+f' g(V,Y')+f\,g(V',Y')+f\,g(RV,Y)\right]\;\mathrm ds=\\
=\int_0^t&\left[-f\,g(V',Y')-f\,g(V,Y'')+f\,g(V',Y')+f\,g(V,RY)\right]\;\mathrm ds=0.
\end{split}
\end{equation}
Similarly, if $\hat V\in\hat{\mathcal K}_0$ and $f\cdot Y(0)\in \hat{\mathcal S}_0$ are fixed,
$f(0)=f(1)=0$, since $g(\hat V',Y(0))$ is constant, from \eqref{eq:defC0} we have:
\[C_0(\hat V,f\cdot Y(0))=\int_0^1f'g(\hat V',Y(0))\;\mathrm du=0,\]
which concludes the proof.
\end{proof}

\begin{corKer}\label{thm:corKer}
For all $t>0$, the kernel of the restriction of $I_t$ to $\mathcal K_t$ equals the kernel of
$I_t$ in $\mathcal H_t$ (see formula \eqref{eq:KerIt}); moreover, $C_0$ is non degenerate in
$\hat{\mathcal K}_0$.
\end{corKer}
\begin{proof}
Let $t\in\,]0,1]$ be fixed. From~\eqref{eq:symmetry}, \eqref{eq:KerIt} and 
\eqref{eq:defcalKt} 
it follows immediately that $\mathrm{Ker}(I_t)\subset\mathcal K_t$, hence
$\mathrm{Ker}(I_t)\subset \mathrm{Ker}(I_t\vert_{\mathcal K_t})$.

For the opposite inclusion, observe that, if $V\in \mathrm{Ker}(I_t\vert_{\mathcal K_t})$,
then $I_t(V,W)=0$ for all $W\in \mathcal K_t$, and, by Lemma~\ref{thm:orthogonal},
also $I_t(V,W)=0$ for all $W\in\mathcal S_t$. By Corollary~\ref{thm:corsoma} it then
follows that $I_t(V,W)=0$ for all $W\in\mathcal H_t$, proving that 
$\mathrm{Ker}(I_t)\supset \mathrm{Ker}(I_t\vert_{\mathcal K_t})$.

Similarly, $\mathrm{Ker}(C_0)=\mathrm{Ker}(C_0\vert_{\hat{\mathcal K}_0})$.
Since $g$ is non degenerate, from \eqref{eq:defC0} it is easy to see that
$C_0$ is non degenerate in $\mathcal H$, which proves that $C_0$ is
non degenerate in $\hat{\mathcal K}_0$.
\end{proof}
We now look at the representation of the bilinear forms
$I_t$ and $C_0$ as self-adjoint operators. 
We start with the following general observation. 

If $B:H^1([0,1],\R^n)\times H^1([0,1],\R^n)\mapsto\R$
is a bilinear form obtained by the restriction of a continuous bilinear form on
$C^0([0,1],\R^n)\times C^0([0,1],\R^n)$, then, since
the inclusion $H^1\mapsto C^0$ is compact, it follows that 
$B$ is represented by a {\em compact operator\/}
on $H^1([0,1],\R^n)$.
 
We can now prove the following:
\begin{Idcomp}\label{thm:Idcomp}
For all $t>0$, $I_t$ is represented by a self-adjoint bounded
linear operator on $\mathcal K_t$ which is of the form $L+K$, where
$L$ is a {\em positive\/} isomorphism of $\mathcal K_t$ and $K$ is compact. 
Also, the restriction of $C_0$ to
$\hat{\mathcal K}_0$ is represented by a compact perturbation of the 
identity map of  $\hat{\mathcal
K}_0$.
\end{Idcomp}
\begin{proof}
Let $t\in\,]0,1]$ be fixed; from \eqref{eq:defIt},
\eqref{eq:ggr} and \eqref{eq:defcalKt} we write
$I_t$ on $\mathcal K_t$ as follows:
\begin{equation}\label{eq:ItonKt}
\begin{split}
I_t(V,W)&=\int_0^t \gr_s(V'(s),W'(s))\;\mathrm ds+\\&+2\int_0^t\frac{[C_V+g(V(s),Y'(s))][
C_W+g(W(s),Y'(s))]}{g(Y(s),Y(s))} \;\mathrm ds+\\
&+\int_0^tg(R(s)[V(s)],W(s))\;\mathrm ds
- g(S[V(0)],W(0)).
\end{split}
\end{equation}
Now, the bilinear form on $\mathcal K_t$ given by the first integral in
\eqref{eq:ItonKt} is a Hilbert space inner product on $\mathcal K_t$,
and therefore it is represented by the identity operator on $\mathcal K_t$.

We now observe that the bounded linear operator \[V\mapsto C_V=\frac1t\int_0^t[g(V',Y)-g(V,Y')]
\,\mathrm ds\] from $H^1([0,t],\R^n)$ to $\R$ has a continuous extension
to $C^0([0,1],\R^n)$. Namely:
\[C_V=\frac 1t\int_0^t\left[g(V',Y)-g(V,Y')\right]\;\mathrm ds=\frac1t\left[
g(V,Y)\big\vert_{0}^t-2\int_0^tg(V,Y')\;\mathrm ds\right],\]
and the latter expression is clearly continuous with respect to the uniform topology.
It follows that the bilinear form on $\mathcal K_t$ given by the second integral
of formula \eqref{eq:ItonKt} has a continuous extension to $C^0([0,1],\R^n)$,
and we have observed that this implies that it is represented by a compact
operator on $\mathcal K_t$. The terms in the last line of formula \eqref{eq:ItonKt} are also
continuous in the $C^0$-topology, and again the corresponding bilinear form
 is represented by a compact operator on $\mathcal K_t$, which proves the first part of 
the Proposition.

As to the bilinear form $C_0$ on $\hat{\mathcal K}_0$, observe that,
by definition of $\hat{\mathcal K}_0$ (see formula \eqref{eq:defcalK0}), 
if $\hat V\in \hat{\mathcal K}_0$
then the quantity $\gr_0(\hat V',Y(0))=-g(\hat V',Y(0))$ is constant, and thus:
\begin{equation}\label{eq:richiamo}
\gr_0(\hat V',Y(0))=\int_0^1\gr_0(\hat V',Y(0))\;\mathrm du=-\gr_0(\hat V(0),Y(0)).
\end{equation}
Then, for $\hat V,\hat W\in \hat{\mathcal K}_0$, it is:
\begin{equation}\label{eq:serve}
\begin{split}
C_0(\hat V,\hat W)&=\int_0^1\gr_0(\hat V'(u),\hat W'(u))\;\mathrm du+\\&-2\,
\frac{\gr_0(\hat V(0),Y(0))\,\gr_0(\hat W(0),Y(0))}{\gr_0(Y(0),Y(0))}.
\end{split}
\end{equation}
Again, the integral in the above formula is a Hilbert space inner product
in $\hat{\mathcal K}_0$, and the last term is continuous in the $C^0$-topology,
which proves that $C_0$ is represented by a compact perturbation of a positive
isomorphism of $\hat{\mathcal K}_0$.
\end{proof}

\begin{initial}\label{thm:initial}
The index of $C_0$ in $\hat{\mathcal K}_0$ is equal to the index of
the restriction of $g$ to the subspace $P$:
\begin{equation}\label{eq:initial}
n_-(C_0\vert_{\hat{\mathcal K}_0})=n_-(g\vert_P).
\end{equation}
\end{initial}
\begin{proof}
Let $P=P_+\oplus P_-$ be a direct sum decomposition of $P$, with $g\vert_{P_+}$ positive
definite and $g\vert_{P_-}$ negative definite (recall that $g$ is non degenerate on $P$).
Then, it is easy to see that 
we have a direct sum decomposition $\hat{\mathcal K}_0=\hat{\mathcal K}_+\oplus
\hat{\mathcal K}_-$, where:
\begin{equation}\label{eq:defK+}
\hat{\mathcal K}_+=\Big\{\hat V\in \hat{\mathcal K}_0:\hat V(0)\in P_+\Big\},
\end{equation}
and
\begin{equation}\label{eq:defK-}
\hat{\mathcal K}_-=\Big\{\hat V:[0,1]\mapsto\R^n\ \text{affine function}\;\big\vert\; \hat
V(0)\in P_-,\ \hat V(1)=0\Big\}.
\end{equation}
Clearly, $\mathrm{dim}(\hat{\mathcal K}_-)=\mathrm{dim}(P_-)=n_-(g\vert_{P})$;
to conclude the proof, it suffices to show that $C_0$ is positive semi-definite
on $\hat{\mathcal K}_+$ and negative definite in $\hat{\mathcal K}_-$.

If $\hat V\in \hat{\mathcal K}_-$, $\hat V\ne0$,
then $\hat V(u)=v_0(u-1)$ for some $v_0\in P_-$, $v_0\ne0$, and for all $u\in[0,1]$;
then, from \eqref{eq:defC0}, we have:
\[C_0(\hat V,\hat V)=\int_0^1g(\hat V',\hat V')\;\mathrm du=g(v_0,v_0)<0.\]
If $\hat V\in\hat{\mathcal K}_+$, then, by \eqref{eq:serve}, we have:
\begin{equation}\label{eq:C0inK+}
C_0(\hat V,\hat V)=\int_0^1 \gr_0(\hat V',\hat
V')\;\mathrm du-2\,\frac{\gr_0(\hat V(0),Y(0))^2}{\gr_0(Y(0),Y(0))}.
\end{equation}
Since $\hat V(1)=0$ ad the function $v\mapsto\gr_0(v,v)$ is convex
in $\R^n$, we use the Jensen's
inequality  to prove the following:
\begin{equation}\label{eq:immSc}
\gr_0(\hat V(0),\hat V(0))=\gr_0(\int_0^1 \hat V'\;\mathrm du,\int_0^1 \hat V'\;\mathrm du ) \le
\int_0^1\gr_0(\hat V'(u),\hat V'(u))\;\mathrm du.
\end{equation}
Finally, from \eqref{eq:C0inK+} and \eqref{eq:immSc}  we obtain:
\[
C_0(\hat V,\hat V)\ge \gr_0(\hat V(0),\hat V(0))-2\,\frac{%
\gr_0(\hat V(0),Y(0))^2}{\gr_0(Y(0),Y(0))}=g(\hat V(0),\hat V(0))\ge0,
\]
which concludes the proof.
\end{proof}
We summarize the above results in the next theorem:
\begin{summarize}\label{thm:summarize}
Let $g$ be a nondegenerate symmetric bilinear form on $\R^n$ with $n_-(g)=1$,
$R:[0,1]\mapsto\mathcal L(\R^n)$ be a $C^1$-map of $g$-symmetric linear operators
on $\R^n$, $P$ a $g$-nondegenerate subspace of $\R^n$ and $S:P\mapsto P$ be a $g$-symmetric
linear map on $P$. Suppose that the differential equation $V''=RV$ admits a solution
$Y$ satisfying $g(Y,Y)<0$ on $[0,1]$. Let $\mathcal K$ be the subspace of $H^1_P([0,1],\R^n)$
consisting of those $V$ such that $g(V',Y)-g(V,Y')$ is constant on $[0,1]$; assume
that $g$ is non degenerate on each $\mathbb J[t]$. Then
\begin{equation}\label{eq:finalformula}
n_-(I_1\vert_{\mathcal K})=n_-(g\vert_P)+\sum_{t\in\,]0,1[}\sgn(t)-n_-(g\vert_{\mathbb J[1]^\perp}),
\end{equation}
where the objects $I_1$ and $\mathbb J[t]$ are defined in \eqref{eq:defIt} and \eqref{eq:defbbJt}.
\qed
\end{summarize}
\end{section}
\begin{section}{On the Nondegeneracy Assumption. The Maslov Index.} 
\label{sec:nondeg}
In this section we will discuss the nondegeneracy assumption for the restriction
of the bilinear form $g$ on the spaces $\mathbb J[t]$ defined in \eqref{eq:defbbJt}, 
and which is essential for the proof of Theorem~\ref{thm:indexthm}.

As we have observed, this assumption guarantees that the set of $(P,S)$-focal
instants is discrete (Proposition~\ref{thm:discrete}); however, it is important
to observe that, even when the number of $(P,S)$-focal instants is finite, such
assumption cannot be removed from the statement of Theorem~\ref{thm:indexthm}
(see \cite[Section~7]{MPT}). 

A natural substitute for the term $\sum_{t\in\,]0,1[}\sgn(t)$ appearing in formula
\eqref{eq:indextheorem} in the case that $g$ is possibly degenerate on some
$\mathbb J[t]$ is the so called {\em Maslov index\/} of the differential problem
\eqref{eq:MS} and \eqref{eq:IC}, denoted by $\maslov(g,R,P,S)$ (see \cite{Hel, MPT} for
details). The Maslov index $\maslov(g,R,P,S)$ is defined whenever $t_0=1$ is not a $(P,S)$-focal instant.
It is  an integer number computed as the intersection number of a continuous curve with
a subvariety  of codimension one of the Lagrangian Grassmannian of a symplectic space.

For the reader's convenience, we sketch briefly the formal definition of
$\maslov$; the proofs and further details on our approach may be found,
in \cite{MPT}.  
Consider the differential problem in $\R^n$  given by \eqref{eq:MS} and \eqref{eq:IC}.
Using the bilinear form $g$, one considers the symplectic form $\omega$ in $\R^{2n}$ 
given by:
\[\omega((v_1,v_2),(w_1,w_2))=g(v_1,w_2)-g(v_2,w_1).\]
It is an easy observation that, if $V$ and $W$ are solutions of \eqref{eq:MS},
then the quantity
$\omega((V(t),V'(t)),(W(t),W'(t))$ is constant in $[a,b]$; moreover,
if $V$ and $W$ are in $\mathbb J$, then this constant is null 
(see formula~\eqref{eq:symmetry}). A subspace $L$ of $\R^{2n}$ is said
to be {\em isotropic\/} with respect to $\omega$ if $\omega$ is null
on $L\times L$; the  space \[L=\Big\{(v_1,v_2)\in\R^{2n}:v_1\in P,\ v_2+S[v_1]\in P^\perp\Big\}\]
is a {\em Lagrangian subspace\/} of the symplectic 
space $(\R^{2n},\omega)$, which is a {\em maximal\/} isotropic subspace of $\R^{2n}$ 
(necessarily $n$-dimensional). The set $\Lambda$ consisting of all
the Lagrangian subspaces of the symplectic space $(\R^{2n},\omega)$ is
a compact, connected, analytic embedded submanifold of the Grassmannian $G_n(\R^{2n})$,
called the {\em Lagrangian Grassmannian\/} of  $(\R^{2n},\omega)$. 

By what has been observed, for all $t\in[a,b]$, the subspace of $\R^{2n}$ given by:
\[L(t)=\Big\{(V(t),V'(t)):V\in\mathbb J\Big\}\]
is Lagrangian, hence the differential problem \eqref{eq:MS} and \eqref{eq:IC}
defines a continuous curve in $\Lambda$. Considering the Lagrangian subspace
of $\R^{2n}$:
\[L_0=\{0\}\oplus\R^n,\]
it is an easy observation that an instant $t_0\in\,]a,b]$ is $P$-focal if and only if
$L(t_0)\cap L_0\ne\{0\}$, i.e., if and only if $L(t_0)$ and $L_0$ are {\em transversal}. 
One then considers the subset $\Lambda_0\subset\Lambda$ consisting of those
Lagrangians that are transversal to $L_0$; $\Lambda_0$ is a dense open subset
of $\Lambda$ which is contractible. The first relative homology group
with integer coefficients $H_1(\Lambda,\Lambda_0;\Z)$ is computed in \cite{MPT} as:
\[H_1(\Lambda,\Lambda_0;\Z)\simeq\Z.\]
The continuous curve $L(t)$ in $\Lambda$ defined by our differential
problem does not define a homology class in $H_1(\Lambda,\Lambda_0;\Z)$, because
its initial point is never in $\Lambda_0$; moreover, its final point is
in $\Lambda_0$ precisely when $t_0=b$ is not a $P$-focal point.
Let's assume that $t_0=b$ is not a $P$-focal point;
by Proposition~\ref{thm:discrete}, if we consider the restriction
$L_\varepsilon$ of the curve $L(t)$ to an interval of the form $[a+\varepsilon,b]$, 
with $\varepsilon>0$ small enough, then we have a well defined continuous curve in $\Lambda$
with endpoints in $\Lambda_0$. The relative homology class of this curve is easily
seen not to depend on the choice of the small $\varepsilon$. The
Maslov index $\maslov(g,R,P,S)$ is defined to be the relative homology
class of $L_\varepsilon$ in $H_1(\Lambda,\Lambda_0;\Z)$.

Such index equals the sum $\sum_{t\in\,]0,1[}\sgn(t)$ when the non degeneracy assumption
for $g$ is satisfied (\cite[Theorem~5.1.2]{MPT}). Moreover, the essential property of
$\maslov$ is that, since it is a topological invariant, it is {\em stable\/} by $C^0$-small
perturbations of the data $(g,R,P,S)$ (\cite[Theorem~5.2.1]{MPT}). As an immediate
application of the uniform stability of $\maslov$, we obtain immediately that the result
of Theorem~\ref{thm:summarize} can be extended to the case that $R$ is only continuous,
provided that the instant $t_0=1$ is not $(P,S)$-focal, by replacing the term $\sum_{t\in\,]0,1[}\sgn(t)$
in \eqref{eq:finalformula} with the Maslov index $\maslov(g,R,P,S)$.

Using a similar perturbation argument, we now want to push the result of Theorem~\ref{thm:summarize}
beyond the assumption of non degeneracy for $g$. To this aim, we argue as follows.

Let's assume that a set of data $(g,R,P,S)$ is given in $\R^n$, with $n_-(g)=1$, and suppose that the
following assumptions are satisfied:
\begin{itemize}
\item[(a)] $g$ is non degenerate on $P$;
\item[(b)] the differential equation $V''\!\!=\!RV$ admits a solution $Y$ satisfying \mbox{$g(Y,Y)<0$}
in $[0,1]$;
\item[(c)] the instant $t_0=1$ is not $(P,S)$-focal.
\end{itemize}
If $g'$ is a symmetric bilinear form on $\R^n$ which is sufficiently close to
$g$ and $P'$ is a subspace of $\R^n$ sufficiently close to $P$ (in the sense of the Grassmannian
of subspaces of $\R^n$), 
then clearly $n_-(g')=1$ and $g'$ is non degenerate on $P'$. So, the assumption (a) above
is stable by small perturbations.

Moreover, standard results on the continuous dependence from the data for ordinary differential
equations guarantee that also the assumptions (b) and (c) above are stable by
uniformly small perturbations of the objects $g$, $R$, $P$ and $S$.

Finally, to complete the argument, we need to prove that it is possible to 
produce arbitrarily $C^0$-small perturbations of the data $(g,R,P,S)$ for
which the restriction of $g$ to the spaces $\mathbb J[t]$ is non degenerate
for all $t\in\,]0,1]$. It is easy to prove that such perturbations of the 
Morse--Sturm problem \eqref{eq:MS} and \eqref{eq:IC}
exist in the more general class of {\em linearized Hamiltonian systems},
where some of the results of this paper and of \cite{MPT} still hold 
in a more general form.
In this class, the set of systems for which the non degeneracy
assumption is $C^0$-dense. Since both the Morse index and the Maslov
index are stable by uniformly small perturbations (see~\cite{MPT}), we
obtain the following extension of Theorem~\ref{thm:summarize}:
\begin{extension}\label{thm:extension}
Let $(g,R,P,S)$ be a set of data for the Morse--Sturm problem \eqref{eq:MS} and \eqref{eq:IC}.
Suppose that the following assumptions are satisfied:
\begin{itemize}
\item $n_-(g)=1$;
\item $R$ is continuous;
\item $t_0=1$ is not a $(P,S)$-focal instant;
\item the equation $V''=RV$ admits a solution $Y$ satisfying $g(Y,Y)<0$ on $[0,1]$;
\end{itemize}
Let $\mathcal K$ be the subspace of $H^1_P([0,1],\R^n)$
consisting of those $V$'s such that the quantity $g(V',Y)-g(V,Y')$ is constant a.e.\ on 
$[0,1]$. Then
\begin{equation}\label{eq:finalformulabis}
n_-(I_1\vert_{\mathcal K})=n_-(g\vert_P)+\maslov,
\end{equation}
where $\maslov=\maslov(g,R,P,S)$ is the Maslov index of the Morse--Sturm problem
and   $I_1$ is the bilinear form on $H^1_P([0,1],\R^n)$ defined in \eqref{eq:defIt}.
\qed
\end{extension}
\end{section}

\begin{section}{The Lorentzian Morse Index Theorem}
\label{sec:Morse}
The main motivation for studying extensions of the Morse--Sturm theory in
the case of non positive metrics $g$ comes from the applications to the
geodesic problem in semi-Riemannian geometry. In this section we discuss
the case of Lorentzian manifolds, and in particular we show how Theorem~\ref{thm:extension}
can be interpreted as a generalization of the classical Morse Index Theorem.

We introduce the following geometrical setup.

Let's assume that $(\M,g)$ is a Lorentzian manifold, $n=\mathrm{dim}(\M)$,  and that
$\gamma:[0,1]\mapsto\M$ is a geodesic, i.e., $\nabla_{\dot\gamma}\dot\gamma=0$, where $\nabla$ is the
covariant derivative of the Levi--Civita connection of $g$. We denote by $R$ the curvature tensor
of $\nabla$, chosen with the following sign convention:  
$R(X,Y)=\nabla_X\nabla_Y-\nabla_Y\nabla_X-\nabla_{[X,Y]}$.  

Let $\mathcal P$ be a smooth submanifold of $\M$, with $\gamma(0)\in\mathcal P$, 
$\dot\gamma(0)\in T_{\gamma(0)}\mathcal P^\perp$, $\gamma(1)=q$, and  assume that $g$ is non degenerate
on $T_{\gamma(0)}\mathcal P$; we say that $\mathcal P$ is non degenerate at $\gamma(0)$. 
The {\em second fundamental form\/} of
$\mathcal P$ at 
$\gamma(0)$ in the direction $n$ is the symmetric bilinear form ${S}_n:
T_{\gamma(0)}{\mathcal P}\times T_{\gamma(0)}{\mathcal P}\mapsto\R$ given by:
\[{S}_n(v_1,v_2)={g}(\nabla_{v_1}V_2,n),\]
where $V_2$ is any extension of $v_2$ to a vector field on $\mathcal P$.
Since $g$ is non degenerate on $T_{\gamma(0)}\mathcal P$, then there exists 
a linear operator, still denoted by $S_n$, on $T_{\gamma(0)}{\mathcal P}$,
such that ${S}_n(v_1,v_2)=g({S}_n[v_1],v_2)$
for all $v_1,v_2\in T_{\gamma(0)}{\mathcal P}$.

A {\em Jacobi field\/} along $\gamma$ is a smooth vector field $J$ along $\gamma$
that satisfies the Jacobi equation 
\begin{equation}\label{eq:Jacobieq}
\nabla_{\dot\gamma}^2J+R(\dot\gamma,J)\,\dot\gamma=0;\end{equation}
a $\mathcal P$-Jacobi field is a Jacobi field $J$ along $\gamma$ that satisfies the 
initial conditions:
\begin{equation}\label{eq:ICgeom}
J(0)\in T_{\gamma(0)}\mathcal P,\quad \big[\nabla_{\dot\gamma(0)}J+S_{\dot\gamma(0)}[J(0)]\big]\in
T_{\gamma(0)}\mathcal P^\perp.
\end{equation}
The {\em index form\/} $I_{\{\gamma,\mathcal P\}}$  is the symmetric bilinear form 
defined on the vector space $\mathcal H_{\{\gamma,\mathcal P\}}$ consisting of those
piecewise smooth vector fields $V$ along $\gamma$ such that $V(0)\in T_{\gamma(0)}\mathcal P$
and $V(1)=0$, defined by:
\begin{equation}\label{eq:indexformgeom}
\begin{split}
I_{\{\gamma,\mathcal P\}}(V,W)=&\;\int_0^1\!\Big[g(\nabla_{\dot\gamma}V,\nabla_{\dot\gamma}W)+
g(R(\dot\gamma,V)\,\dot\gamma,W)\Big]\;\mathrm dt+\\&-  g(S_{\dot\gamma(0)}[V(0)],W(0)).
\end{split}
\end{equation}
It is easy to see that a vector field $V\in\mathcal H_{\{\gamma,\mathcal P\}}$ is a $\mathcal P$-Jacobi
field if and only if it is in the kernel of $I_{\{\gamma,\mathcal P\}}$. A point $\gamma(t_0)$
is said to be a $\mathcal P$-focal point along $\gamma$ if there exists a non zero
$\mathcal P$-Jacobi field along $\gamma$ vanishing at $t_0$; the multiplicity of a $\mathcal P$-focal
point is the dimension of the vector space of all $\mathcal P$-Jacobi fields along $\gamma$
vanishing at $t_0$. 
If the initial submanifold $\mathcal P$ reduces to a fixed point of $\M$, in which case
the $\mathcal P$-Jacobi fields along $\gamma$ are simply the Jacobi fields
vanishing at $t=0$, then the focal points are also called {\em conjugate points}.
If $\gamma$ is either timelike or lightlike, in which case
$\mathcal P$ is necessarily a spacelike submanifold of $\M$ at $\gamma(0)$, 
then there are only a finite number of $\mathcal P$-focal points along $\gamma$, and
their number, with multiplicity, is defined to be the {\em geometric index\/}
of the geodesic $\gamma$ (see~\cite{PT}). 

The geodesic $\gamma$ is a critical point of the action functional:
\begin{equation}\label{eq:action}
f(z)=\frac12\int_0^1g(\dot z,\dot z)\;\mathrm dt,
\end{equation}
defined on the set $\Omega_{\{\mathcal P,q\}}$
of paths $z:[0,1]\mapsto\M$ such that $z(0)\in \mathcal P$ and
$z(1)=q$; the space $\mathcal H_{\{\gamma,\mathcal P\}}$ can be seen as the 
tangent space of $\Omega_{\{\mathcal P,q\}}$ at $\gamma$
and the bilinear form $I_{\{\gamma,\mathcal P\}}$ is the {\em second variation\/} of $f$
at $\gamma$. Hence, the index of $I_{\{\gamma,\mathcal P\}}$ in $\mathcal H_{\{\gamma,\mathcal P\}}$
is the Morse index of the functional  $f$ at the critical point  $\gamma$; moreover, $\gamma$
is a non degenerate critical point of $f$ precisely when the point $q$ is not
$\mathcal P$-focal along $\gamma$. 
\smallskip

The Morse index of $f$ at $\gamma$ is not finite, due to the indefiniteness of the
metric $g$. However, the theory developed in the previous sections indicate that
we can determine a finite index carrying some geometric information about $\gamma$
provided that we restrict the bilinear form $I_{\{\gamma,\mathcal P\}}$ to a suitable
subspace of $\mathcal H_{\{\gamma,\mathcal P\}}$.
\smallskip

To describe how the geometrical problem fits into the theory of Morse--Sturm systems
discussed in the previous sections, we consider a {\em trivialization\/}
of the tangent bundle $T\M$ along $\gamma$ by means of a  family
$\{E_1,\ldots,E_n\}$ of parallel vector fields along $\gamma$.

The map $V=\sum_i\lambda_i\cdot E_i\mapsto (\lambda_1,\ldots,\lambda_n)$ gives
an isomorphism of $\mathcal H_{\{\gamma,\mathcal P\}}$ with the vector
space of all piecewise smooth $\R^n$-valued functions on $[0,1]$. Since each
$E_i$ is parallel, the covariant derivative of vector fields along $\gamma$
correspond to the usual differentiation in $\R^n$; moreover, the Lorentzian metric
$g$ is carried to a constant nondegenerate bilinear form on $\R^n$, still denoted
by $g$, with $n_-(g)=1$. For each $t\in[0,1]$, the map  \[\R^n\simeq T_{\gamma(t)}\M\ni v\mapsto
R(\dot\gamma(t),v)\,\dot\gamma(t)\in T_{\gamma(t)}\M\simeq\R^n\] is 
given by a $g$-symmetric linear operator on $\R^n$, still denoted by $R(t)$.
Finally, the tangent space $T_{\gamma(0)}\mathcal P$ corresponds to a $g$-nondegenerate
subspace $P$ of $\R^n$, and the second fundamental form $S_{\dot\gamma(0)}$ gives a
$g$-symmetric linear map $S:P\mapsto P$.  

The bilinear form $I_{\{\gamma,\mathcal P\}}$ is carried into the bilinear form
$I_1$, defined in the set of piecewise smooth $\R^n$-valued functions on $[0,1]$,
given by formula \eqref{eq:defIt}. Since $I_1$ has a continuous extension
to the Hilbert space $H^1_P([0,1],\R^n)$, an easy density argument shows
that the index of $I_{\{\gamma,\mathcal P\}}$ on $\mathcal H_{\{\gamma,\mathcal P\}}$
is equal to the index of $I_1$ on $H^1_P([0,1],\R^n)$. 
The Jacobi equation \eqref{eq:Jacobieq} becomes the Morse--Sturm system \eqref{eq:MS},
the initial conditions \eqref{eq:ICgeom} are read into \eqref{eq:IC}, and we have
translated our Lorentzian geodesic problem into the Morse--Sturm problem \eqref{eq:MS}
and \eqref{eq:IC}.

Clearly, the space $\mathbb J$\/ defined in \eqref{eq:defbbJ} corresponds to the set
$\mathcal J_{\mathcal P}$ of $\mathcal P$-Jacobi fields, and the $(P,S)$-focal instants are precisely
the $\mathcal P$-focal points along $\gamma$. The space $\mathcal J_{\mathcal P}[t]\subset
T_{\gamma(t)}\M$ is defined to be the set of values at $t$ of the fields in $\mathcal J_{\mathcal P}$;
the {\em signature\/} $\sgn(\gamma(t_0))$ of the $\mathcal P$-focal point 
$\gamma(t_0)$ is defined to be the signature of the metric $g$ on the space 
$\mathcal J_{\mathcal P}[t_0]^\perp$; the $\mathcal P$-focal point $\gamma(t_0)$
is said to be {\em positive, null or negative\/} according to whether $\sgn(\gamma(t_0))$
is positive, null or negative.\footnote{The reader should observe that we are using
a  terminology  slightly different from the one adopted in \cite{Hel},
where it is defined a {\em timelike, a null and a spacelike index\/} for each conjugate
point.}

The important observation here
is that, if $\gamma$ is causal, i.e., timelike or lightlike, then the restriction
of the metric $g$ to the space $\mathcal J_{\mathcal P}[t]^\perp$ is always positive definite,
so that the signature of a $\mathcal P$-focal point coincides with its
multiplicity. 
This implies in particular that the Maslov index of $\gamma$ coincides precisely with
the geometrical index of $\gamma$.
\smallskip

Under the assumption that the point $\gamma(1)$ is not $\mathcal P$-focal along
$\gamma$, we can therefore apply Theorem~\ref{thm:extension} to the geometrical
problem, obtaining the following generalization of the Morse Index Theorem
for Lorentzian geodesics with variable initial endpoint:
\begin{MorseLorentzian}\label{thm:MorseLorentzian}
Let $(\M,g)$ be a Lorentzian manifold, $\mathcal P\subset\mathcal M$ a smooth
submanifold, $\gamma:[0,1]\mapsto\M$ a geodesic with $\gamma(0)\in\mathcal P$ and $\dot\gamma(0)\in
T_{\gamma(0)}\mathcal P^\perp$. Assume that the following  hypotheses are satisfied:
\begin{itemize}
\item there exists a timelike Jacobi field $Y$ along $\gamma$;
\item $\mathcal P$ is non degenerate at $\gamma(0)$;
\item $\gamma(1)$ is not $\mathcal P$-focal along $\gamma$.
\end{itemize}
Then, denoting by $\mathcal K^\gamma$ the space of (piecewise smooth) vector fields
$V$ along $\gamma$ satisfying $V(0)\in T_{\gamma(0)}\mathcal P$, $V(1)=0$ and
$g(\nabla_{\dot\gamma}V,Y)-g(V,\nabla_{\dot\gamma}Y)\equiv C_V$ (constant), the index
of $I_{\{\gamma,\mathcal P\}}$ on $\mathcal K^\gamma$ is finite, and the following equality
holds:
\begin{equation}\label{eq:LorentzianMorse}
n_-(I_{\{\gamma,\mathcal P\}}\vert_{\mathcal K^\gamma})=n_-(g\vert_{T_{\gamma(0)}\mathcal P})+\maslov(\gamma).
\end{equation}
Moreover, if $\gamma$ is causal, then $\maslov(\gamma)$ equals the geometric index
of $\gamma$.\qed
\end{MorseLorentzian}
Observe that the quantity on the right hand side of \eqref{eq:LorentzianMorse} does {\em not\/}
depend
on the choice of the timelike Jacobi field $Y$, hence the index of $I_{\{\gamma,\mathcal P\}}$
on the space $\mathcal K^\gamma$ is also independent on the choice of $Y$.
We also remark that, if $\gamma$ is a timelike geodesic, then one can take as a timelike
Jacobi field $Y$ the tangent field $\dot\gamma$. It is easy to see that, in this case,
the space $\mathcal K^\gamma$ consists precisely of those vector fields
along $\gamma$ that are pointwise orthogonal to $\dot\gamma$. Hence, Theorem~\ref{thm:MorseLorentzian}
gives a generalization of the Timelike Morse Index Theorem of \cite[Theorem~10.27]{BEE}.
\smallskip

An important class of examples where the assumption on the existence of a timelike
Jacobi field  along any geodesic is satisfied is given by the {\em stationary Lorentzian
manifolds}, i.e., Lorentzian manifolds admitting a timelike Killing vector field. In this case,
a timelike Jacobi vector field along every geodesic
is given by the restriction of any timelike Killing field (see~\cite[Lemma~9.26, p.\ 252]{ON}).

It is interesting to observe that, for non positive definite metrics, as we can
deduce from equation \eqref{eq:LorentzianMorse}, the Morse
index of the action functional at a given geodesic $\gamma$ may be strictly positive
even in the case that $\gamma$ has no focal points. This happens precisely when the
initial submanifold $\mathcal P$ is non spacelike. For a better understanding of
this fact, one can consider the following simple but instructive example.

\begin{exsimple}\label{thm:exsimple}
Let $(\M,g)$ be the two-dimensional flat Minkowski space, with metric $\mathrm dx^2-\mathrm dy^2$.
Let $\gamma(t)=(t,0)$, $t\in[0,1]$, and  let $\mathcal P$ denote the one-dimensional
timelike submanifold of $\M$ given by the $y$-axis; we are in the situation described in the
hypotheses of Theorem~\ref{thm:MorseLorentzian}, considering $Y=\frac\partial{\partial y}$
as the timelike Jacobi field along $\gamma$. 
Clearly, there are no $\mathcal P$-focal points along $\gamma$, and both the
curvature tensor $R$ of $g$ and the second fundamental form $S$ of $\mathcal P$ are null.

We have $n_-(g\vert_{T_{(0,0)}\mathcal P})=1$; 
the space $\mathcal K^\gamma$ consists of vector fields of the form $V=a(t)\frac\partial{\partial x}
+b(t)\frac\partial{\partial y}$, with $a(0)=a(1)=b(1)=0$ and $b'(t)\equiv C_V$ constant on $[0,1]$.
For  $V\in\mathcal K^\gamma$, the value of the index form $I_{\{\gamma,\mathcal P\}}(V,V)$
is computed easily as:
\[I_{\{\gamma,\mathcal P\}}(V,V)=\int_0^1\big[a'(t)^2-b'(t)^2\big]\;\mathrm dt=
\int_0^1 a'(t)^2 \;\mathrm dt-C_V^2.\]
If we consider $a\equiv 0$, we get a one-dimensional subspace of $\mathcal K^\gamma$
on which $I_{\{\gamma,\mathcal P\}}$ is negative definite; on the other hand,
if we consider $b\equiv 0$ and $a$ arbitrary, we get a complementary 
subspace where $I_{\{\gamma,\mathcal
P\}}$ is positive definite, thus $n_-(I_{\{\gamma,\mathcal P\}}
\vert_{\mathcal K^\gamma})\!=~\!\!1$.
\end{exsimple}
It is fairly easy to give examples of $\mathcal P$-focal points in stationary
Lorentzian manifolds of every causal type.  Examples of positive focal or conjugate
points are easily constructed by considering causal geodesics, or spacelike geodesics
admitting a parallel timelike Jacobi field along them
(see~Example~\ref{thm:exparallel} below).
In the next example we construct elementary examples of negative and null
focal points in manifolds with flat metric.
\begin{excausal}\label{thm:excausal}
Consider the Minkowski plane $\R^2$ endowed with the flat metric
$g=\mathrm dx^2-\mathrm dy^2$; let $\gamma(t)=(t,0)$ be the (spacelike) geodesic
segment on the $x$-axis, and let $\mathcal P$ denote the parabola
 through the origin given by the equation $y^2+2x=0$. Then, $\gamma$ is orthogonal to
$\mathcal P$ at $(0,0)=\gamma(0)$; the second fundamental form of $\mathcal P$
at $(0,0)$ is easily computed as 
\[S_{\dot\gamma(0)}\left(\frac\partial{\partial y}\right)=\frac\partial{\partial y},\]
so that $J(t)=(t-1)\frac\partial{\partial y}$ is a $\mathcal P$-Jacobi field along
$\gamma$ which vanishes at $t=1$. Clearly, $\gamma(1)=(1,0)$ is a $\mathcal P$-focal
point of multiplicity one along $\gamma$, and $\sgn(\gamma(1))=
\sgn(g\vert_{\R\cdot J'(1)})=\sgn(g\vert_{\R\cdot \frac\partial{\partial
y}})=-1$.

To construct an example of a null $\mathcal P$-focal point we now consider
the three-dimensional flat Minkowski space $\M=\R^3$ with metric $g=\mathrm dx^2+\mathrm
dy^2-\mathrm dz^2$ and the spacelike geodesic  $\gamma(t)=(t,0,0)$, $t\in[0,1]$.
Let $\mathcal P$ be any smooth surface through the origin such that the tangent
plane $T_{(0,0,0)}\mathcal P$ is the $yz$-plane and such that the second fundamental form
$S_{\dot\gamma(0)}$ of $\mathcal P$ at $(0,0,0)$ satisfies%
\footnote{of course, such submanifold $\mathcal P$ exists; see for instance
\cite[Lemma~2.3.2]{MPT} for details on how to construct a smooth submanifold
of a semi-Riemannian manifold when its tangent space and its second fundamental form
is assigned at one point.}
\[S_{\dot\gamma(0)}(\frac\partial{\partial y}+\frac\partial{\partial z})=
\frac\partial{\partial y}+\frac\partial{\partial z}.\]
Arguing as before
it is easy to verify that $J(t)=(t-1)(\frac\partial{\partial y}+\frac\partial{\partial z})$
is a $\mathcal P$-Jacobi field along $\gamma$, $J(1)=0$, $J'(1)$ is the lightlike
vector $\frac\partial{\partial y}+\frac\partial{\partial z}$, and $\gamma(1)$ is
a null $\mathcal P$-focal point along $\gamma$.
\end{excausal}
\begin{remsemi}\label{thm:remsemi}
Theorem~\ref{thm:MorseLorentzian} can be extended to the
case of geodesics in semi-Rie\-mann\-ian manifolds $(\M,g)$, with
$g$ of arbitrary index $n_-(g)=k\ge1$. In this case, given a geodesic
$\gamma$ in $\M$, one needs
to assume the existence of $k$   Jacobi fields $J_1,\ldots, J_k$ 
along $\gamma$ generating a $k$-dimensional timelike distribution 
along $\gamma$, and satisfying the  relations
$g(\nabla_{\dot\gamma}Y_i,Y_j)-g(Y_i,\nabla_{\dot\gamma}Y_j)=0$ for all
$i,j=1,\ldots,k$. One considers the space $\mathcal K^\gamma$ of vector fields
$V$ along $\gamma$ satisfying $V(0)\in T_{\gamma(0)}\mathcal P$,
$V(1)=0$ and $g(\nabla_{\dot\gamma}V,Y_i)-g(V,\nabla_{\dot\gamma}Y_i)\equiv C_V^{(i)}\ \text{(const.)}$
for all $i=1,\ldots,k$. Then, if $\gamma(1)$ is not $\mathcal P$-focal,
the index of $I_{\{\gamma,\mathcal P\}}$ on $\mathcal K^\gamma$ equals $\maslov(\gamma)+
n_-(g\vert_{T_{\gamma(0)}\mathcal P})$.

Examples of semi-Riemannian manifolds where the theory applies
are given by those manifolds admitting a family of Killing vector fields
$Y_1,\ldots,Y_k$ that generate a $k$-dimensional timelike distribution on $\M$, 
and satisfying the commutation relations $[Y_i,Y_j]=0$ for all $i,j=1,\ldots,k$.
A variational theory for geodesics in this kind of manifolds is presented in
\cite{GPS}. All the results presented in this paper can be extended to
this more general situation.
\end{remsemi}
We now discuss the case of conjugate points along Lorentzian geodesics satisfying
the hypotheses of Theorem~\ref{thm:MorseLorentzian}, and so we assume that
the initial manifold $\mathcal P$ reduces to a single point.
This means that the $\mathcal P$-Jacobi fields along $\gamma$ are simply 
the Jacobi fields vanishing at $t=0$. We denote by $I_\gamma$ the index
form along $\gamma$ relative to the choice of a trivial initial
manifold.

The first observation is that, in this situation, if $Y$ is parallel
along the geodesic $\gamma$, then the conjugate points along $\gamma$ are
isolated, and they are all positive.
\begin{exparallel}\label{thm:exparallel}
Suppose that $Y$ is a parallel timelike Jacobi field along the geodesic
$\gamma$; this means that $\nabla_{\dot\gamma}Y=0$, and so $\nabla_{\dot\gamma}^2Y=
R(\dot\gamma,Y)\,\dot\gamma=0$. 

If $J$ is Jacobi, then $g(\nabla_{\dot\gamma}J,Y)-g(J,\nabla_{\dot\gamma}Y)=
g(\nabla_{\dot\gamma}J,Y)$ is constant on $[0,1]$, hence $\frac{\mathrm d^2}{\mathrm dt^2}
g(J,Y)=g(\nabla_{\dot\gamma}^2J,Y)=0$, and $g(J,Y)$ is an affine function on $[0,1]$.

If $\gamma(t_0)$ is conjugate to $\gamma(0)$ along $\gamma$, and $J$ is a non trivial
Jacobi field along $\gamma$ vanishing at $0$ and $t_0$, then it must be $g(J,Y)\equiv0$,
and so $g(\nabla_{\dot\gamma}J,Y)\equiv0$. It is $\mathcal J_{\mathcal P}[t_0]^\perp=
\{\nabla_{\dot\gamma}J(t_0):J\ \text{Jacobi, with }\ J(0)=J(t_0)=0\}$, and it
follows that $\mathcal J_{\mathcal P}[t_0]^\perp\subset Y(\gamma(t_0))^\perp$. Since
$Y$ is timelike, it follows that the restriction of the metric $g$ to 
$\mathcal J_{\mathcal P}[t_0]^\perp$ is positive definite, which implies that
the conjugate point $\gamma(t_0)$ is  isolated and that its
signature $\sgn(\gamma(t_0))$ is equal to its multiplicity. 
Hence, the Maslov index of $\gamma$ coincides
with its geometric index. In this case, Theorem~\ref{thm:MorseLorentzian} tells us that, 
if $\gamma(1)$
is not conjugate to $\gamma(0)$ along $\gamma$, the index of $I_{\gamma}$ on 
$\mathcal K^\gamma$ is equal to the geometric index of $\gamma$.
\end{exparallel}

Let's assume now that the geodesic $\gamma$ satisfies the assumptions of
Theorem~\ref{thm:MorseLorentzian} and the non degeneracy assumption of 
Proposition~\ref{thm:discrete}. It is an easy observation that there
cannot be {\em too many\/} negative conjugate 
points along $\gamma$. For example, if $\gamma(t_0)$ is a negative
conjugate point, i.e., $\sgn(\gamma(t_0))=-1$, then the Maslov index
$\maslov(\gamma\vert_{[0,t_0-\varepsilon]})$ must be strictly positive 
for $\varepsilon>0$ small enough. This follows immediately from the fact 
that, by Theorem~\ref{thm:MorseLorentzian}, if $\varepsilon>0$ is small 
enough, it must be 
\[n_-(I_\gamma\vert_{\mathcal
K^\gamma_{t_0+\varepsilon}})=\maslov(\gamma\vert_{[0,t_0+\varepsilon]})=\maslov(\gamma\vert_{[0,t_0-\varepsilon]})-1
\ge0.\]
In particular, the first conjugate point along $\gamma$ is {\em never\/} negative.
\smallskip

If $\mathrm{dim}(\M)=2$, then the metric $-g$ is Lorentzian in $\M$. This simple
observation allows to get some interesting consequences, like the following:
\begin{dim2}\label{thm:dim2}
Let $(\M,g)$ be a two dimensional Lorentzian manifold and let $\gamma:[0,1]\mapsto\M$ be
a spacelike geodesic in $\M$. Suppose that there exists a timelike Jacobi field
along $\gamma$. Then, there are no conjugate points along $\gamma$.
\end{dim2}
\begin{proof}
The curve $\gamma$ is clearly a timelike geodesic in the opposite Lorentzian
manifold $(\M,-g)$ with the same conjugate points. We know that all the conjugate
points along a causal geodesic are positive, hence $\gamma$ has only negative
conjugate points in $(\M,g)$. Then, there cannot be any conjugate point, because
the sum of their signatures must be non negative integer.
\end{proof}
By the same argument, it is easy to see that if $\gamma$ is a spacelike geodesic
in a two-dimensional Lorentzian manifold $\mathcal M$, starting orthogonally to a 
one-dimensional (necessarily timelike) submanifold $\mathcal P$ of $\mathcal M$,
then there is {\em at the most one\/} $\mathcal P$-focal point along $\gamma$,
which must be negative (see Example~\ref{thm:excausal}).
It is well known that conjugate points cannot occur along lightlike geodesics
in two-dimensional Lorentzian manifolds (see~\cite{BEE}). However,
we remark that spacelike (or timelike) geodesics in two-dimensional Lorentzian
Lorentzian manifolds may have conjugate points. For instance, in the conformally flat
metric $e^{t^2}(\mathrm dx^2-\mathrm dt^2)$ on $\R^2$, the curve $\gamma(\tau)=(\tau,0)$
is a spacelike geodesic, and the Jacobi equation along $\gamma$ is for
the vector field $J=(v,w)$ is given by the system
\[v''=0,\quad w''+w=0.\]
Clearly, the point $\gamma(\pi)$ is conjugate to $\gamma(0)$ along $\gamma$.
\smallskip
 
We leave unanswered the following questions:
\begin{enumerate}
\item do there exist examples of (spacelike) Lorentzian geodesics satisfying the hypotheses
of Theorem~\ref{thm:MorseLorentzian} for which the set of $\mathcal P$-focal (or conjugate)
points is not discrete?
\item can a (spacelike) Lorentzian geodesic satisfying the hypotheses of 
Theorem~\ref{thm:MorseLorentzian} really have one negative conjugate point?
\item suppose that $\gamma$ is a (spacelike) geodesic satisfying the hypotheses of 
Theorem~\ref{thm:MorseLorentzian} and having one or more isolated conjugate
point for which the non degeneracy assumption of Proposition~\ref{thm:discrete}
is not satisfied; is it still true that the Maslov index of $\gamma$ is given by the
sum of the signatures of its conjugate points?
\end{enumerate}
If one does not require the assumptions of Theorem~\ref{thm:MorseLorentzian} all
the above questions have easy answers (see~\cite{MPT}): the first two questions
have a positive answer and the third one has a negative answer.
\smallskip

We conclude this section with the remark that a Lorentzian version of
the Morse Index Theorem for the two variable endpoints (see \cite{Kal} for the Riemannian
case) can be easily deduced from Theorem~\ref{thm:MorseLorentzian}. When the final
endpoint of $\gamma$ is allowed to vary on a submanifold $\mathcal Q$ of $\mathcal M$,
the index of the second variation of the action functional at $\gamma$ is given by
the sum of the right hand side of equation \eqref{eq:LorentzianMorse} and a term
that measures the {\em relative convexity\/} of $\mathcal Q$ with respect to $\mathcal P$.
The details are found in \cite[Theorem~2.7, Remark~2.10]{PT}.

\end{section}

\begin{section}{The Global Morse Relations for Geodesics in Stationary Lorentzian
Manifolds}\label{sec:globalMorse}
In this section we want to develop an infinite dimensional Morse
theory for the geodesics joining two fixed points $p$ and $q$ 
in a stationary Lorentzian manifold $(\M,g)$, in the spirit of \cite{M}
and using the modern terminology of \cite{Bott}.
The main goal of this theory is to give estimates on the number of geodesics 
having a given index; these estimates are given in terms of the topology of 
the space of (continuous)  curves  joining $p$ and $q$ in $\M$. 
The basic reference for most of the material discussed in this section is 
\cite{GP}; we will make full use of the results proven in that article.

As customary, if $I\subseteq\R$ is any interval, we will denote by $H^1(I,\R^n)$
the Sobolev space of absolutely continuous curves  $z:I\longmapsto\R^n$  
such that the integral
$\int_I\vert\dot z\vert^2\;{\rm d}t$ is finite, 
where   $\vert\cdot\vert$ denotes the
Euclidean norm in $\R^n$. 
\smallskip

Given {\em any\/} differentiable manifold $N$, the set $H^1([0,1],N)$ is
defined as the set of all absolutely continuous curves $z:[0,1]\longmapsto N$
such that, for every local chart $(V,\varphi)$ on $N$, with $\varphi:U\longmapsto \R^n$
a diffeomorphism, and for every closed subinterval
$I\subseteq[0,1]$ such that $z(I)\subset V$, it is $\varphi\circ z\in H^1(I,\R^n)$. 
For all differentiable manifold $N$, with $\mathrm{dim}(N)=n$, the set $H^1([0,1],N)$ 
has the structure
of an infinite dimensional manifold, modeled on the Hilbert space $H^1([0,1],\R^n)$.
We will denote by $TN$ the tangent bundle of $N$ and by $\pi:TN\mapsto N$ the
canonical projection; for $p\in N$, $T_pN=\pi^{-1}(p)$ denotes the tangent
space of $N$ at $p$.
A vector field along a curve $z:[0,1]\mapsto N$ is a map $\zeta:[0,1]\mapsto TN$ with
$\pi(\zeta(t))=z(t)$ for all $t$. Given any $z\in H^1([0,1],N)$, the tangent
space $T_zH^1([0,1],N)$ is identified with the set:
\[T_zH^1([0,1],N)=\Big\{\zeta\in H^1([0,1],TN):\zeta\ \text{vector field along}\ z\Big\},\]
which is an infinite dimensional vector space, with a topology that makes it into
a {\em Hilbertable\/} space.

Let's assume that $(\M,g)$ is a Lorentzian manifold which admits 
a timelike Killing vector field, denoted by $Y$. We assume that
$Y$ is complete; let $p$ and $q$ be fixed points in $\M$.
We introduce the following space:
\[\Omega_{p,q}=\Big\{z\in H^1([0,1],\M):z(0)=p,\ z(1)=q\Big\},\]
It is well
known that  $\Omega_{p,q}$ has the structure of an infinite
dimensional Hilbertian submanifold of $H^1([0,1],\M)$; the action functional $f$, defined in
\eqref{eq:action}, is smooth on $\Omega_{p,q}$ and its critical points are precisely the geodesics
in $\M$ between $p$ and $q$. We say that $p$ and $q$ are non conjugate in $\M$ if
they are not conjugate along every geodesic in $\M$ joining them.

For all geodesic $\gamma$ in $\M$ we have a conservation law
$g(\dot\gamma,Y)\equiv c_\gamma\,\text{(constant)}$. 
Now, if we consider the subset $\Omega_{p,q}^Y$ of $\Omega_{p,q}$
consisting of curves $z$ satisfying $g(\dot z,Y)\equiv\text{const.}$, then clearly
the geodesics in $\Omega_{p,q}$ belong to $\Omega_{p,q}^Y$.
It is proven in \cite{GP} that $\Omega_{p,q}^Y$ is a smooth submanifold of
$\Omega_{p,q}$, and that  $f$ has the same critical
points in $\Omega_{p,q}$ and in $\Omega_{p,q}^Y$.

By differentiating the expression $g(\dot z,Y)=\text{const.}$ with respect to $z$, using the Killing
property of $Y$ it is easy to see that the tangent space $T_z\Omega_{p,q}^Y$ is given by
the Hilbert space of $H^1$-vector fields along $z$ satisfying $V(0)=V(1)=0$ and such that
the quantity $g(\nabla_{\dot z}V,Y)-g(V,\nabla_{\dot z}Y)$ is constant a.e.\ on $[0,1]$.

Hence, if $\gamma$ is a critical point for $f$ in $\Omega_{p,q}^Y$, i.e., a geodesic
between $p$ and $q$, the tangent
space $T_\gamma\Omega_{p,q}^Y$ is a completion of the space $\mathcal K^\gamma$
of Theorem~\ref{thm:MorseLorentzian}, and the index of $I_{\gamma}$ in
$\mathcal K^\gamma$ is equal to the Morse index of the functional $f$ at the
critical point $\gamma\in\Omega_{p,q}^Y$.

Such index can therefore be interpreted as the number of {\em essentially different directions\/}
in which $\gamma$ can be deformed, in the class of curves $z$ joining $p$ with $q$
and satisfying $g(\dot z,Y)=\text{const.}$, in order to obtain a curve with smaller action.

Let $\mathcal C^1_{p,q}$ denote the following space:
\[\begin{split}\mathcal C^1_{p,q}=\Big\{z:[0,1]\mapsto&\M\ \text{piecewise}\ C^1 
:\\ & z(0)=p,\ z(1)=q,\ g(\dot z,Y)\equiv c_z\
\text{(constant)}\Big\};
\end{split}\]
we give the following completeness condition for the sublevels 
 of the restricted action functional.
\begin{defcprec}\label{thm:defcprec}
Given $c\in\R$, we say that $\mathcal C_{p,q}^1$ is $c$-precompact if every sequence
$\{z_n\}_{n\in\N}\subset\mathcal C^1_{p,q}$ such that $f(z_n)\le c$ has a uniformly
convergent subsequence.
\end{defcprec}
The $c$-precompactness property, which is given intrinsically in Definition~\ref{thm:defcprec}, 
can be studied by means of suitable bounds of the metric coefficients with respect to the 
coordinates of a given atlas on $\M$. A wide class of examples of stationary Lorentzian 
manifolds $(\M,g)$ for which the $c$-precompactness assumption is satisfied by all choices of
$p$, $q$ and $c$ is given in \cite{GP}. We emphasize that the $c$-precompactness for
stationary Lorentzian manifold plays the role of the completeness assumption in
Riemannian geometry; for this and other analogies with the classical Riemannian
theory we refer to \cite{GP}, where it is also discussed the relation between the 
$c$-precompactness and the property of {\em global hyperbolicity}. 

We recall that, given a topological space $X$, an algebraic field $\corpo$ and a natural 
number $i$, the $i$-th {\em Betti number\/} $\beta_i(X;\corpo)$ of $X$ relative to $\corpo$
is the $\corpo$-dimension of the $i$-th singular homology vector space
$H_i(X;\corpo)$ of $X$ with coefficients in $\corpo$.  The {\em Poincar\'e polynomial\/}
$\mathfrak P_\lambda(X;\corpo)$ of $X$ with coefficients in $\corpo$ is the formal
power series in $\lambda\in\corpo$ given by:
\begin{equation}\label{eq:defpoincare}
\mathfrak P_\lambda(X;\corpo)=\sum_{i=0}^\infty\beta_i(X;\corpo)\,\lambda^i.
\end{equation}
The global Morse relations provide relations between the set of all the
geodesics joining $p$ and $q$ in $\M$ with the topology of the space of all
continuous curves joining $p$ and $q$ in $\M$, given in terms of the Betti numbers
and the Poincar\'e polynomial of this space.  
A key point for the infinite dimensional Morse theory is the so called
{\em Palais--Smale\/} condition. We recall that a smooth functional $F$
on a manifold $X$ endowed with a Finsler structure
is said to satisfy the Palais--Smale condition at the
level $c\in\R$ if every sequence $\{x_n\}_{n\in\N}\subset X$ such that:
\begin{itemize}
\item[(a)] $\lim\limits_{n\to\infty}F(x_n)= c$;
\item[(b)] $\lim\limits_{n\to\infty} \Vert{\mathrm d}F(x_n)\Vert=0$,
\end{itemize}
has a converging subsequence in $X$. 

The $c$-precompactness condition given in Definition~\ref{thm:defcprec}
is the key assumption for the proof of the global Morse relations,
which are given in the following

\begin{globalMorse}\label{thm:globalMorse}
Let $(\M,g)$ be a Lorentzian manifold. Suppose that $\M$ admits a 
{\em complete\/} timelike Killing vector field $Y$, and assume that $p$ and
$q$ are two points of $\M$ such that the following hypotheses are satisfied:
\begin{itemize}
\item $p$ and $q$ are not conjugate in $\M$;
\item $\mathcal C^1_{p,q}$ is $c$-precompact for all $c\in\R$.
\end{itemize}
Let $\Omega_{p,q}^0$ denote the space of all continuous curves $z:[0,1]\mapsto\M$ joining
$p$ and $q$ in $\M$, endowed with the topology of uniform convergence, and let
$\mathcal G_{p,q}$ denote the set of all geodesics in $\M$ between $p$ and $q$.
Then, for all field $\corpo$ there exists a formal power series $Q_\corpo(\lambda)$ 
in the variable $\lambda$, with coefficients in $\N\bigcup\{+\infty\}$ such that the
following identity between formal power series is satisfied:
\begin{equation}\label{eq:powerseries}
\sum_{z\in\mathcal G_{p,q}}\lambda^{\maslov(z)}=\poincare_\lambda(\Omega_{p,q}^0;\corpo)+(1+\lambda)\,
Q_\corpo(\lambda).
\end{equation}
\end{globalMorse}
\begin{proof}
Let $f_Y$ denote the restriction of the action functional
$f$ to the manifold $\Omega_{p,q}^Y$; as we have observed,
$\Omega_{p,q}^Y$ is a smooth submanifold of $\Omega_{p,q}$
and the critical points of $f_Y$ on $\Omega_{p,q}^Y$ are
precisely the geodesics joining $p$ and $q$ in $\M$.

We endow $\Omega_{p,q}^Y$ with the following Riemannian structure.
We consider an auxiliary Riemannian metric $\gr$ on $\M$, and
for all $z\in \Omega_{p,q}^Y$ we define a Hilbert space inner
product $\rip\cdot\cdot$ in $T_z\Omega_{p,q}^Y$ by:
\begin{equation}\label{eq:riemstructure}
\rip VV=\int_0^1\gr({\nabla_{\dot z}V},{\nabla_{\dot z}V})\;\mathrm dt.
\end{equation}

Using the $c$-precompactness assumption, as well as the density of
$\mathcal C^1_{p,q}$ in $\Omega_{p,q}^Y$, the following facts are
proven in \cite{GP}:
\begin{enumerate}
\item $f_Y$ is {\em bounded from below}, i.e., there exists $D\in\R$
such that $f(z)\ge D$ for all $z\in\Omega_{p,q}^Y$;
\item for all $c\in\R$, the sublevel $f_Y^c=\big\{z\in\Omega_{p,q}^Y:f(z)\le c\big\}$
is a {\em complete\/} metric subspace of $\Omega_{p,q}^Y$;
\item for all $c\in\R$, $f_Y$ satisfies the Palais--Smale condition at the level $c$
when $\Omega_{p,q}^Y$ is endowed with the Finsler structure given by \eqref{eq:riemstructure}.
\end{enumerate}
Finally, the condition that $p$ and $q$ be non conjugate in $\M$ implies
that $f_Y$ is a {\em Morse functional}, i.e., all its critical points in $\Omega_{p,q}^Y$
are non degenerate. Namely, as we have already observed, the second variation
of $f_Y$ at any geodesic $\gamma$  is given by the restriction of the
index form $I_\gamma$, and its kernel in $T_\gamma\Omega_{p,q}^Y$ coincides with
the set of Jacobi fields along $\gamma$ vanishing at the endpoints. If $p$ and $q$
are non conjugate in $\M$, then $I_\gamma$ has trivial kernel, and $f_Y$ is a 
Morse functional.

\noindent\ \ 
Then, by standard results of Global Analysis on Manifolds (see for instance~\cite{MW}),
denoting by $m(z,f_Y)$ the Morse index of the critical point $z$ of $f_Y$,
we have the following Morse relations. For all field $\corpo$ there exists a 
formal power series $Q_\corpo(\lambda)$  in the variable $\lambda$, 
with coefficients in $\N\bigcup\{+\infty\}$ such that the following identity 
between formal power series is satisfied:
\begin{equation}\label{eq:MorseRelParziali}
\sum_{z\in\mathcal G_{p,q}}\lambda^{m(z,f_Y)}=\mathfrak P_\lambda(\Omega^Y_{p,q};\corpo)+
(1+\lambda)\,Q_\corpo(\lambda).
\end{equation}
By Theorem~\ref{thm:MorseLorentzian}, for all $z\in\mathcal G_{p,q}$
we have $m(z,f_Y)=\maslov(z)$; moreover, since $Y$ is complete, it is proven
in \cite{GP} that the spaces $\Omega_{p,q}$ and $\Omega_{p,q}^Y$ are
homotopically equivalent, which implies that $\mathfrak P_\lambda(\Omega^Y_{p,q};\corpo)=
\mathfrak P_\lambda(\Omega_{p,q};\corpo)$ for all field $\corpo$. Finally,
also the spaces $\Omega_{p,q}$ and $\Omega_{p,q}^0$ have the same homotopy type
(see~\cite{M}), and so the Morse relations \eqref{eq:powerseries} are easily obtained from
\eqref{eq:MorseRelParziali}.
\end{proof}
\end{section}



\end{document}